\documentclass[12pt]{amsart}
\usepackage{amsmath,amssymb,amscd,array, mathrsfs }

\usepackage{amsmath,amscd,amssymb,amsthm,array}
\usepackage[all,cmtip]{xy} 
\usepackage{color}

\usepackage{amsmath,amssymb,mathrsfs,amsthm, tikz-cd,mathrsfs}
\usepackage{comment}
\def\url#1{\expandafter\s

\tring\csname #1\endcsname}

\def\mmat #1,#2,#3,#4,{\text{\small\arraycolsep=3pt $
\begin{pmatrix}#1&#2\\#3&#4\end{pmatrix}$}}

\usepackage{hyperref}
\usepackage{capt-of}

\usepackage{multirow}

\usepackage{comments}

\def\mmat #1,#2,#3,#4,{\text{\small\arraycolsep=3pt $
\begin{pmatrix}#1&#2\\#3&#4\end{pmatrix}$}}

\usepackage{lscape}
\usepackage{tikz-cd}
\usepackage{enumerate}

\usepackage{DLdef1}
\usepackage{multicol}

\def\mmat #1,#2,#3,#4,{\text{\small\arraycolsep=3pt $
\begin{pmatrix}#1&#2\\#3&#4\end{pmatrix}$}}

\renewcommand {\ssbegin}[2][*]
 {\refstepcounter{subsection}%
\if#1*
\addcontentsline{toc}{subsection}{\thesubsection.\hskip 1pc #2}%
\else
\addcontentsline{toc}{subsection}{\thesubsection.\hskip 1pc #2. #1}%
\fi
 \def \secno {\gdef \secno {}{\ssecfont
\thesubsection.\hskip 2ex}%
 }%
 \begin{#2}}

\renewcommand {\sssbegin}[2][*]
 {\refstepcounter{subsubsection}
\if#1*
\addcontentsline{toc}{subsubsection}{\thesubsubsection.\hskip 1pc #2}%
\else
\addcontentsline{toc}{subsubsection}{\thesubsubsection.\hskip 1pc #2. #1}
\fi
 \def \secno {\gdef \secno {}{\ssecfont \thesubsubsection.\hskip 2ex}%
 }%
 \begin{#2}}

\renewcommand {\parbegin}[2][*]
 {\refstepcounter{paragraph}
\if#1*
\addcontentsline{toc}{paragraph}{\theparagraph.\hskip 1pc #2}%
\else
\addcontentsline{toc}{paragraph}{\theparagraph.\hskip 1pc #2. #1}
\fi
 \def \secno {\gdef \secno {}{\ssecfont \theparagraph.\hskip 2ex}%
 }%
 \begin{#2}}

\newcommand {\ce}{{\text{CE}}}

\newcommand{\N}{\mathbb{N}^{*}}

\rmnameii{etr}{etr}
\rmnameii{evv}{ev}

\DeclareMathOperator{\K}{\mathbb{K}}

\newcommand{\blue}{\color{blue}}

\newcommand{\black}{\color{black}}
\newcommand{\Z}{\mathbb{Z}}
    
\begin{document}

\title{On the Classification of Nilpotent Lie Superalgebras over Fields of Characteristic 2}
\author{Muhammad Uzair Haroon}
\begin{abstract}
The salient result of this paper is the  complete classification of all  nilpotent Lie superalgebras in dimensions up to four over an algebraically closed field of characteristic 2. Our method adapts the Skjelbred-Sund method to show that every $n$-dimensional nilpotent Lie superalgebra arises as a one-dimensional central extension of an $(n-1)$-dimensional one, and  determines precisely when two such extensions yield isomorphic superalgebras. Carrying out this construction for $\dim\le4$, we obtain a finite, non‐redundant list of pairwise non-isomorphic examples.
\end{abstract}
\maketitle


\thispagestyle{empty}
\setcounter{tocdepth}{2}
\tableofcontents

\section{Introduction} \label{Introduction}
\vspace{0.25 em}

\subsection{Lie Superalgebras in Characteristic 2.}\label{1} Lie superalgebras arose in the 1970s as $\mathbb{Z}/2\mathbb{Z}$-graded analogues of Lie algebras, motivated by supersymmetry in physics. In this setting, elements are assigned degrees, which leads to the appearance of signs in the defining relations. Roughly speaking, a Lie superalgebra is a $\Z/2\Z$-graded vector space $L$ equipped with a parity-preserving bilinear map $[\cdot,\cdot]:L\times L\to L$ called the \textit{Lie bracket}, which is super anti-symmetric and satisfies a graded version of the Jacobi identity. Although \textit{simple} Lie superalgebras have been classified by Kac (\cite{K}), the classification of other classes of Lie superalgebras, such as nilpotent ones is limited to a few particular cases, see Section \ref{1.2}.

Over a field of characteristic 2, a Lie superalgebra is a $\Z/2\Z$-graded vector space $L=L_{\ev}\oplus L_{\od}$ such that the even part $L_{\ev}$ is a Lie algebra, the odd part $L_{\od}$ is a two-sided  $L_{\ev}$-module due to symmetry, and there is a map $s:L_{\od}\to L_{\ev}$, called the \textit{squaring} which satisfies certain consistency relations, including a modified Jacobi identity (see Section \ref{deflsa}). The first formal instance of this definition, appears in \cite{LeD}, however, earlier analogues may be found in \cite{BMY} and \cite{NR} with the name ``strongly anti-commutative graded algebras".

Due to the introduction of the squaring map, characteristic 2 is a novel case that requires novel techniques. In fact, certain Lie superalgebras in characteristic 2 have no equivalents in other characteristics (see \cite{BGL1, BGLLS, BLLSq, LeD}), the standard example being the four mutually non-isomorphic versions of the Hamiltonian Lie superalgebra. Thus, the classification of both simple and nilpotent cases is an open problem. However, a number of breakthroughs have permeated the setting. Lie superalgebras in characteristic 2 that admit a Cartan matrix have been classified in \cite{BGL1}, albeit under the assumption that they each accede a Dynkin diagram with exactly one odd node. Furthermore, the classification of simple Lie superalgebras in characteristic 2 has been reduced to the classification of simple Lie algebras in characteristic 2 (by either queerification or ``method 2"), as described in \cite{BLLSq,BGLLS}. The latter classification is itself yet open and a difficult problem in its own right (see \cite{S1,S2}). 

A full cohomology theory was developed and introduced in \cite{BM} to investigate the structural properties of these superalgebras, although preliminary descriptions may be found in \cite{BGLL}. This cohomology is different from the usual Chevalley-Eilenberg complex for Lie (super)algebras over fields of characteristic 0 and similarly does not have any analogues in characteristic $p>2$. We only require $1$ and $2-$cocycles for the purposes of this paper, and their constructions are given in Section \ref{cohom}. As an application, Lagrangian extensions of Lie superalgebras in characteristic 2 have been studied and classified in dimension 4 (see \cite{BBE}) utilizing this cohomology theory.

\subsection{Classification of Nilpotent Lie (super)algebras.} \label{1.2} The classification of \emph{simple} Lie algebras has been achieved by Kac \cite{K} who generalized the classical results of Cartan and Killing. Over a field of characteristic $p>0$ the problem is still open except for a few particular cases, as mentioned in Section \ref{1}. In the case where the (super)algebra is nilpotent, the problem is more challenging and a full classification remains an open question. Due to nilpotence however, it is possible to obtain an explicit classification by detailing the relations on the basis elements.

The classification of nilpotent Lie algebras over any field has been completed
for a long time up to dimension 5. For characteristic 0, the dimension 6 classification was achieved by Morozov in 1958 (\cite{M}), whereas the result was achieved for any arbitrary field in 2012, see \cite{CGS,DG}. On the other hand, the classification of nilpotent Lie superalgebras is known only in some particular cases: for any arbitrary field in characteristic $0$, those with a (super)multiplier rank $\leq 2$ have been classified (see \cite{LZ} for both the definition and the classification result), and for $\K=\mathbb{R}$ or $\mathbb{C}$, all nilpotent Lie superalgebras have been classified up to dimension $5$, as found in \cite{H}. Over a \emph{perfect} field of characteristic $p>0$, the \emph{restricted} $p$-nilpotent Lie algebras were classified up to dimension 4 by Schneider and Usefi in \cite{SU} (including the $p=2$ case), followed by dimension 5 for a perfect field of characteristic $p>3$ by Darijani and Usefi in \cite{DU}. The latter result, however, was partially incorrect, and the full list for $p>3$ was found by Maletesta and Siciliano in \cite{MS}. 

Over any arbitrary field of characteristic $p\geq3$, the restricted $p$-nilpotent Lie \emph{superalgebras} were classified up to dimension 4 by Bouarroudj and Ehret in \cite{BE}. The method they used built on the classical Skjelbred-Sund method \cite{SS} that can be described as follows. Given a Lie superalgebra $L$ of dimension $n-1$, they compute all \textit{central extensions} of $L$ by a one-dimensional vector space $V$. That is, they define the Lie bracket and the $p$-map on the direct sum $L \oplus V$ in such a way that $V$ lies in the center of the extended superalgebra. The obstructions for the extended bracket to satisfy the Jacobi identity are controlled by the cohomology of $L$ with coefficients in the trivial module $V$. This process yields \textit{all} $p$-nilpotent Lie superalgebras of dimension $n$. As a final step, they detect and remove isomorphic duplicates from the resulting list. Due to certain similarities between the squaring map and the $p$-map, we aim to build on their method to classify nilpotent Lie superalgebras in characteristic $p=2$. 

Much work has also been produced recently on the classification of nilpotent non-associative algebras by adapting the Skjelbred-Sund method, including nilpotent Lie triple systems up to dimension 4 in \cite{ABCF1}, nilpotent Lie-Yamaguti algebras up to dimension 4 in \cite{ABCF2}, non-Lie central extensions of 4-dimensional Malcev algebras in \cite{HAC}, and non-associative central extensions of 3-dimensional Jordan algebras in \cite{HA}, amongst others. For a much broader list of algebras that have been classified using the method, see \cite{ABCF1}.

\vspace{1 em}

\subsection{Organization of the paper and main results.} The presentation of the paper, and to a lesser extent the organization, is motivated by that of \cite{BE}. We review basic concepts and fundamental definitions related to Lie superalgebras in characteristic 2 in Section \ref{deflsa}. Section \ref{Engel} is dedicated to proving Engel's Theorem for our setting. We recall the cohomology theory for characteristic 2 in Section \ref{cohom}, detailing particularly only the first and second cohomology groups as relevant. Section \ref{ext} serves as a build-up to Theorem \ref{thm1}, which concerns the relationship of central extensions of Lie superalgebras and the even part of the respective second cohomology group. Moreover, we rectify an oversight from \cite[Theorem 3.5.3]{BE} in Step 5.

The rest of the paper is dedicated to an inductive classification of nilpotent Lie superalgebras. The theory mainly comprises Proposition \ref{prop1} and Theorem \ref{indthm}, which together show that any $n$-dimensional nilpotent Lie superalgebra may be obtained as a central extension of an $(n-1)$-dimensional nilpotent Lie superalgebra by a scalar $2$-cocycle, while determining when such extensions are isomorphic. The computational method for classification is described in Section \ref{method}. Succinctly, we compute the scalar $2$-cocycles of the $2$-dimensional nilpotent Lie superalgebras in \ref{nilclassif2} and determine their equivalence classes with respect to the module structure defined in Section \ref{ModStr}. We proceed to then build the extensions in \ref{4.3} and remove isomorphic copies in \ref{iso3} to yield a \emph{complete, non-redundant} list of nilpotent Lie superalgebras in dimension 3 in Theorem \ref{dim3}. The process is repeated at a much larger scale in Section \ref{class4}, to yield a \emph{complete} list of nilpotent Lie superalgebras in dimension 4 in Theorem \ref{classif4}.

\section{Lie superalgebras in Characteristic 2}

\vspace{0.5 em}

Throughout the text, $\mathbb{K}$ stands for an algebraically closed field of characteristic two. The group of integers modulo $2$ is denoted by $\Z/2\Z$. 

Let $V$ and $W$ be two vector spaces over $\mathbb{K}$. A map $s : V \to W$ is called a \textit{squaring} if
\[ s(\lambda x) = \lambda^2 s(x), \quad \forall \lambda \in \mathbb{K}, \forall x \in V. \]
We also suppose that the following map is bilinear
\[ V \times V \to W \quad (x, y) \mapsto s(x + y) - s(x) - s(y). \]

 Let $V=V_{\bar 0}\oplus V_{\bar 1}$ be a superspace defined over ${\mathbb K}$. The parity of a homogeneous element $v\in V_{\bar{i}}$ is denoted by $|v|:=\bar{i}$. The element $v$ is called \textit{even} if $v\in V_{\bar 0}$ and \textit{odd} if $v\in V_{\bar 1}$. The {\it superdimension} of $V$ is $\mathrm{sdim}(V)=a+b\epsilon$, where $\epsilon^2=1$, and $a=\mathrm{dim}(V_{\bar 0})$, $b=\mathrm{dim}(V_{\bar 1})$. Usually, $\mathrm{sdim}(V)$ is shorthanded as $a\mid b$. The total dimension of $V$, which is $a+b$, is denoted by  $\mathrm{dim}(V)$.  
 
 Throughout the text, all elements are supposed to be homogeneous unless otherwise stated. A linear map $\varphi:V\rightarrow W$ between superspaces is called \textit{even} if $\varphi(V_{\bar{i}})\subset W_{\bar{i}}$ and \textit{odd} if $\varphi(V_{\bar{i}})\subset W_{\overline{i}+\overline{1}}$.

\vspace{0.5 em}

\subsection{Lie superalgebras in characteristic 2.}\label{deflsa}

Following \cite{BMY,LeD}, a \emph{Lie superalgebra} over a field $\K$ of characteristic $p=2$ is a $\Z/2\Z$-graded vector space $L=L_{\ev}\oplus L_{\od}$ such that the even part $L_{\ev}$ is a Lie algebra, the odd part $L_{\od}$ is a $L_{\ev}$-module made two-sided by symmetry (the bracket on $L_{\ev}$ as well as the action of $L_{\ev}$ on $L_{\od}$ are given by the same symbol $[\cdot, \cdot]$), a squaring  $s:L_{\od}\rightarrow L_{\ev}$, such that the bracket of two odd elements is given by:
    \begin{equation}\label{bil}
        [x,y]:=s(x+y)-s(x)-s(y),~~ \forall x,y\in L_{\od}.
    \end{equation} 
The Jacobi identity involving the squaring reads as follows: 
    \begin{equation}\label{jacobsq} 
        [s(x),y]=[x,[x,y]],~\forall x\in L_{\od},~\forall y\in L.    \end{equation}
Such a Lie superalgebra will be denoted by $\left(L,[\cdot,\cdot]_L,s_L\right)$ or by just $L$ if no confusion is possible.

\vspace{0.5 em}

\sssbegin{Remark}Any associative superalgebra $A=A_{\ev}\oplus A_{\od}$ in characteristic $2$ leads to  a Lie superalgebra. The bracket is the usual commutator and the squaring is defined by $s(x) := x\cdot x,~\forall x \in A_\od$. 
\end{Remark}

\vspace{0.5 em}

Let $L$ be a Lie superalgebra and $K\subseteq L$ be a homogeneous linear subspace. We say that 
    \begin{itemize}
        \item $K$ is a \emph{Lie subalgebra} if it is closed under the bracket and under the squaring;
        \item $K$ is an \emph{ideal} if $[k,x]\in K,~\forall k\in K,~\forall x\in L$ as well as $s(k)\in K, ~\forall k\in K_{\od}$. 
    \end{itemize}

Let $L$ and $H$ be two Lie superalgebras. An even linear map $\varphi: L\rightarrow H$ is called a \emph{Lie superalgebras morphism} if
\begin{equation}
        \varphi\left([x,y]_L\right)=[\varphi(x),\varphi(y)]_H,~\forall x,y\in L, \text{ and }
        \varphi\circ s_L(x)=s_H\circ\varphi(x),~\forall x\in L_{\od}. 
\end{equation}

 A \emph{representation} of $L$ in a $\Z/2\Z$-graded vector space $M$ is an even map $\rho: L\rightarrow \text{End}(M)$ satisfying
\begin{eqnarray*}
\rho([x, y])& = &[\rho(x), \rho(y)] \text{ $\forall x, y \in L$; and}\\[1mm]
\rho (s(x)) &=&(\rho(x))^2 \text{ $\forall x \in L_{\od}.$}
\end{eqnarray*}
 Such a $M$ is called a \emph{$L$-module}.\\~

 Let $L=(L,[\cdot,\cdot],s)$ be a Lie superalgebra in characteristic $2$. We denote $$s(L_\od):=\text{Span}\{s(x),x\in L_\od\}.$$

The \emph{center} of $L$ is the ideal
\[
\mathfrak{z}(L) := \Big\{\, x \in L \ \Big|\ [x,y] = 0 \ \text{for all } y \in L, \ \text{and if } x \in L_{\od} \text{ then } s(x)=0 \,\Big\}.
\]

 The \emph{lower central series} of $L$ is defined inductively by
\[
L^{1} := L, \qquad 
L^{k+1} := \text{Span}\{[L, L^{k}] \cup s(L^k_{\od})\} \quad \text{for } k \ge 1.
\]
We say that $L$ is \emph{nilpotent} if there exists an integer $n \ge 1$ such that
\[
L^{n} = 0.
\]

\sssbegin{Proposition}[\text{\cite[Proposition 2.1.3]{BBE}}]\label{classif2}
    Let $L$ be a $2$-dimensional Lie superalgebra over an arbitrary field of characteristic $2$. Then,  $L$ is isomorphic to one of the following superalgebras.
    \end{Proposition}
    
   \begin{itemize} 
     \item \underline{$\sdim(L)=(0|2)$}: $L={\bf L_{0|2}^1}=\left<0|e_1,e_2\right>$.\\
     
    \item \underline{$\sdim(L)=(1|1)$};  $L=\left<e_1|e_2\right>$.
    
    \begin{multicols}{2}

        \begin{enumerate}[$(1)$]
            \item ${\bf L_{1|1}^1} =\left<e_1|e_2;[e_1,e_2]=e_2\right>$;

            \item ${\bf L_{1|1}^2}=\left<e_1|e_2; s(e_2)=e_1 \right>$;

            \columnbreak        
            \item ${\bf L_{1|1}^3} =\left<e_1|e_2\right>$ \textup{(}abelian\textup{)}.
        \end{enumerate}
    \end{multicols}    
    \item \underline{$\sdim(L)=(2|0)$}: $L=\left<e_1, e_2|0\right>$.
        \begin{multicols}{2}
        \begin{enumerate}[$(1)$]
            \item ${\bf L_{2|0}^1} =\left<e_1,e_2|0; [e_1,e_2]=e_2\right>;$

        \columnbreak             
            \item ${\bf L_{2|0}^2} =\left<e_1,e_2|0\right>$ \textup{(}abelian\textup{)}.
           
        \end{enumerate}
        \end{multicols}
    
    \end{itemize} 

\sssbegin{Proposition}\label{nilclassif2}
Let $L$ be a $2$-dimensional nilpotent Lie superalgebra over an arbitrary field of characteristic $2$. Then, $L$ is isomorphic to one of the following Lie superalgebras.
\end{Proposition}

\begin{itemize}
\item \underline{$\mathrm{sdim}(L) = (0|2)$}: 
$L = \mathbf{L^{1}_{0|2}} = \langle 0| e_{1}, e_{2} \rangle \; (\textit{abelian}).$\\

\item \underline{$\mathrm{sdim}(L) = (1|1)$}: 
$L = \langle e_{1}| e_{2} \rangle.$\\
 
\hspace{1em}(1) $\mathbf{L^{2}_{1|1}} = \langle e_{1}|e_{2}; s(e_{2}) = e_{1}\rangle;$  
\hspace{1em}(2) $\mathbf{L^{3}_{1|1}} = \langle e_{1}|e_{2}\rangle \; (\textit{abelian}).$\\

\item \underline{$\mathrm{sdim}(L) = (2|0)$}: 
$L = \langle e_{1}, e_{2}|0\rangle.$\\

\hspace{1em}(1) $\mathbf{L^{2}_{2|0}} = \langle e_{1}, e_{2}|0\rangle \; (\textit{abelian}).$
\end{itemize}

\sssbegin{Remark} For certain examples of simple Lie superalgebras in characteristic 2, see \cite{BGL1, BGLLS, BLLSq}. 
\end{Remark}

\vspace{0.5 em}

\subsection{Engel's Theorem for Lie superalgebras in characteristic \texorpdfstring{$p=2$}{p=2}.}\label{Engel} We prove first Engel's Theorem for this setting. The theorem has long been generalized to Lie superalgebras in characteristic 0 (see \cite{K}), and has also been proven for other classes of algebras and superalgebras (see \cite{AE,AO,HRR,KR,St}). The theorem appears to be recognized among experts for Lie superalgebras in characteristic 2, though one does not find a formal proof for it anywhere. 

Our attempt aims to adapt the standard proof for Lie algebras. In particular, we build on the approach in \cite{F}, though we remain wary and correct some oversights in the original attempt.

\vspace{0.5 em}

\sssbegin{Theorem}
Let $L$ be a Lie superalgebra in characteristic $p=2$. Then, $L$ is nilpotent if and only if all its elements are ad-nilpotent.
\end{Theorem}

It is easy to see why the forward direction holds. Indeed, suppose $L$ is nilpotent. Then, there exists $n \in \N $ such that $L^n=0$. By the definition of nilpotence, this clearly implies, $[x_1,\dots,[x_{n-1},y]\dots]=0$ for all $x_i,y \in L$. Setting all $x_i=x$ implies $(ad_x)^{n-1}(y)=0$ for all $y \in L$. Hence, we have the first implication. For the other direction, we need to introduce some concepts. 

\vspace{0.5 em}

\sssbegin{Definition}
The \textit{normalizer} of a subalgebra $K$ of a Lie superalgebra $L$ is defined as the linear subspace of $L$ 
$$N_{L}(K)=\{x\in L: [x,K] \subseteq K \}.$$
\end{Definition}

\vspace{0.5 em}

It is easy to see that $K \subseteq N_L(K)$. 

In order to prove the reverse implication, we first need to prove three lemmas. 

\vspace{0.5 em}

\sssbegin{Lemma}\label{Lemma3}
$N_L(K) \subseteq L$ is a subalgebra of $L$.
\end{Lemma}

\begin{proof}
Let $K$ be a subalgebra of $L$. We first prove that $N_L(K)$ is homogeneous. Let $x=x_{\ev}+x_{\od} \in N_L(K)$ and $k_{\ev} \in K_{\ev}$, $k_{\od} \in K_{\od}$. 
Then, we have 
\begin{align*}
[x,k_{\ev}]=[x_{\ev},k_{\ev}]+[x_{\od},k_{\ev}]\\
[x,k_{\od}]=[x_{\ev},k_{\od}]+[x_{\od},k_{\od}]
\end{align*}
Since $[x,k_{\ev}] \in K$, we have $[x_{\ev},k_{\ev}] \in K_{\ev}$ and $[x_{\od},k_{\ev}] \in K_{\od}$. Similarly, since $[x,k_{\od}] \in K$, we have $[x_{\ev},k_{\od}] \in K_{\od}$ and $[x_{\od},k_{\od}] \in K_{\ev}$. We hence have that
\begin{align*}
[x_{\ev},k]=[x_{\ev},k_{\ev}]+[x_{\ev},k_{\od}]\in K,\\
[x_{\od},k]=[x_{\od},k_{\ev}]+[x_{\od},k_{\od}]\in K.
\end{align*}
Therefore, $N_L(K)$ is homogeneous. The fact that $N_L(K)$ is closed under the bracket and squaring follows from the Jacobi identity. By taking $x \in K$, $y,z \in N_L(K)$, we have
$$[x,[y,z]]=[y,[x,z]]+[z,[y,x]].$$

Since the left-hand side is in $K$, we have by the right hand  side $[[y,z], K]\subseteq K$, which implies $[y,z] \in N_L(K)$. Moreover, for $k\in K$, $x \in N_L(K)_{\od}$, we have
$$[s(x),k]=[x,[x,k]].$$

Again, using the same logic, we have $s(x)\in N_L(K)$. 
\end{proof}

\vspace{0.5 em}

\sssbegin{Lemma}\label{Lemma 1} Let $V$ be a superspace, and let $L$ be a Lie subalgebra in $\mathfrak{gl}(V)$ of nilpotent endomorphisms. If $V$ is non-zero and finite dimensional, then $V$ must contain a non-zero element $v$ such that $l(v)=0$ for all $l \in L$. 
\end{Lemma}
\begin{proof}
We proceed by induction on $\dim(L)$. Let $\dim(L)=1$. Then, if $l \in L$ is a nilpotent endomorphism, there exists an $n\in \N$ such that $l^{n-1} \neq 0$ but $l^{n}=0$. Choose $v=l^{n-1}(w)\not =0$ for some $w\in V$, then $l(v)=l^n(w)=0$. But $\dim(L)=1$, which implies that for all $l' \in L$, $l'$ is a scalar multiple of $l$, and hence $l'(v)=0$. 

We now take $\dim(L)>1$. We would like to construct an ideal with co-dimension $1$. We start by taking a proper subalgebra $K \subsetneq L$. If $L_{\ev}=0$, then for every $x \in L_{\od}$ we have $s(x)=0$. Hence, $\langle x \rangle = \textrm{Span}(\{x\})$ for any $x\in L_{\od}$ is a proper subalgebra. On the other hand, if $L_{\ev} \neq 0$, then we pick $\langle x \rangle$ for all $x\in L_{\ev}$. Hence, a proper subalgebra exists. 


Since $\dim(K) \lneq \dim(L)$, by the induction hypothesis, there exists a non-zero $x+K \in L/K$ such that $\mathrm{ad}_{k}(x+K)=[k,x]+K=0$ for all $k \in K$. This implies $[k,x] \in K$ for all $k \in K$. Therefore, $x\in N_L(K)$ and $x\not \in K$.  Hence, $K$ is a proper subalgebra of $N_L(K)$. If we let $K$ be the maximal proper subalgebra of $L$, then by this argument and Lemma \ref{Lemma3} we have $N_L(K)=L$. Hence, $K$ is an ideal in $L$. 

Now, we aim to prove that the codimension of $K$ (maximal) has to be equal to $1$ with respect to $L$. We see that if $\dim(L/K)>1$, then as we proved earlier, there exists a one-dimensional subalgebra $H/K \subsetneq L/K$. Then, the preimage of $H/K$ under the canonical map $L \longrightarrow L/K$ is a proper subalgebra of $L$ which contains $K$. But this is a contradiction, as $K$ is maximal. Hence, we have that $\dim(L/K)=1$.

Now, we define $W:= \{v\in V: Kv=0\}$. By assumption, $W$ is non-zero. Moreover, $W$ is stable under the action of $L$. Indeed, let $x \in L$, $y \in K$, $w \in W$, then we have by the commutator
$$yxw=xyw+[x,y]w=0.$$
\ignorespacesafterend
Since $K$ is an ideal and the right term reduces to zero, we have $xw \in W$.

Finally, let $z \in L-K$.  By the definition of $W$ (which is non-empty by assumption) we have that there exists $v\in W$ such that $Kv=0$. Since we know $L/K$ is $1$-dimensional, we just need to prove for a single $z\in L/K$ that it also sends some $v\in W$ to zero. Since every map in $L \subseteq \mathfrak{gl}(V)$ is nilpotent, and $W$ is stable under $L$, if we apply it enough times to any $w \in W$, it will go to zero. We just take $v=z^{n-1}(w)\neq 0$ for some $w\in W$. Then, we know that any map in $L/K$ sends this $v$ to zero as it is 1-dimensional. Hence, we have that any element in $K$ and $L/K$ sends this $v$ to zero, and so it is true for $L$, and we have the claim.
\end{proof}

\sssbegin{Lemma}\label{Lemma2}
If $L/\mathfrak{z}(L)$ is nilpotent, then so is $L$.
\end{Lemma}

\begin{proof}
If $L/\mathfrak{z}(L)$ is nilpotent, then there exists an $n \in \N$ such that $L^n \subseteq \mathfrak{z}(L)$. Then, $L^{n+1}=0$.
\end{proof}

Finally, we prove Engel's Theorem.

\begin{proof}
Suppose $L$ is a Lie superalgebra composed of $\mathrm{ad}$-nilpotent elements. We apply Lemma \ref{Lemma 1} to the Lie superalgebra $\mathrm{ad}(L)$ acting on the space $L$. Then, there exists a non-zero $x\in L$ such that $[L,x]=0$. Next, we prove that $\dim(L/\mathfrak{z}(L))<\dim(L)$.

It is easy to see that such an $x$ may be reduced to the case of being either purely even or purely odd. Indeed, if $x=x_{\ev}+x_{\od}$, then by the bilinearity of the bracket we have for all homogeneous  basis elements $l \in L$
$$[l,x]=[l,x_{\ev}]+[l,x_{\od}]=0.$$
However, since the bracket is an $even$ map, $[l,x_{\ev}] \neq [l,x_{\od}]$, unless both equal 0.

We suppose $x$ is even. Then, such an element is automatically in $\mathfrak{z}(L)$. Now we suppose $x$ is odd. Then we have two cases. Either $s(x)=0$ and therefore $x\in \mathfrak{z}(L)$. Or $s(x)=y$ for some $y \in L_{\ev}$. But then $[y,l]=[x,[x,l]]=0$ for all $l\in L$. Hence, $y \in \mathfrak{z}(L)$. In any case, we see that $\mathfrak{z}(L)$ must contain some element and so $L/\mathfrak{z}(L)$ must have a smaller dimension than $L$.

Using Lemma \ref{Lemma 1} recursively on $L/\mathfrak{z}(L)$, and the case being trivial for $\dim(L)=0$, we see that $L/\mathfrak{z}(L)$ must be nilpotent. By Lemma \ref{Lemma2}, if $L/\mathfrak{z}(L)$ is nilpotent, then so is $L$.
\end{proof}

\vspace{0.5 em}

\subsection{Cohomology Theory of Lie superalgebras in characteristic 2.}\label{cohom}We  recall the cohomology theory  of Lie superalgebras in characteristic 2, particularly only the first and second cohomology groups, since that is all that we will need. The full description of the theory, and in particular that of the higher order cohomology, may be found in \cite{BM}.

\vspace{0.5em}

Let $L$ be a Lie superalgebra in characteristic $2$ and let $M$ be a $L$-module. 

A 1-cocycle on $L$ with values in $M$ is a linear map $\varphi: L \rightarrow M$ such that
\begin{equation}\label{1coc}
\begin{array}{llll}
d_{\ce}^1(\varphi)(x,z)&:=&x\cdot \varphi(z)+z \cdot \varphi(x)+\varphi([x,z])=0, & \forall x,z\in L; \\[2mm] 
\delta^1(\varphi)(x)&:=&x\cdot \varphi(x)+\varphi(s(x))=0,&\forall x\in L_\od.
 \end{array}
\end{equation}
The space of 1-cocycles on $L$ with values in $M$ is denoted by $XZ^1(L;M)$. We also use the notation $\fd^1(\varphi):=(d^1_{\ce}(\varphi),\delta^1(\varphi))$.
\vspace{0.5 em}

A 2-cocycle on $L$ with values in $M$ consists of a pair $(\alpha, \gamma)$  such that: 
\[
\begin{array}{l}
\alpha:L\wedge L \rightarrow M \text{ is a bilinear map}, \text{ and } \gamma:L_{\od} \rightarrow M \text{ is a squaring}
\end{array}
\]
such that $\alpha(x,y)=\gamma(x+y)+\gamma(x)+\gamma(y)$ for all $x,y\in L_{\od}$. Moreover, we require that 
\begin{equation}\label{2coc}
\begin{array}{llll} 
d_{\ce}^2(\alpha)(x,y,z)&:=&\underset{x,y,z}{\circlearrowleft}\bigl(x \cdot \alpha(y,z)+\alpha([x,y],z)\bigl)=0, &\forall x,y,z\in L; \\[2mm] 
\delta^2(\alpha,\gamma)(x,z)&:=&x\cdot \alpha(x,z)+z \cdot \gamma(x)+\alpha(s(x),z)+\alpha([x,z],x)=0,&\forall x\in L_\od,~\forall z\in L,
 \end{array}
\end{equation} 

\vspace{0.2 mm}

The space of 2-cocycles on $L$ with values in $M$ is denoted by $XZ^2(L;M)$. We also use the notation $\mathfrak{d}^2(\alpha,\gamma):=\bigl(d_{\ce}^2(\alpha),\delta^2(\alpha,\gamma)\bigr)$. 

\vspace{0.5 em} 

If $\alpha$ is an even map, then $\mathrm{Im}(\gamma)\subseteq M_{\overline{0}}$, whereas if $\alpha$ is odd, we have $\mathrm{Im}(\gamma) \subseteq M_{\overline{1}}$. Therefore, $XZ^2(L;M)$ can be graded as follows:
\begin{equation}\label{cograd}|(\alpha,\gamma)|:=|\alpha|, \quad \forall (\alpha,\gamma)\in XZ^2(L;M).\end{equation}

\vspace{0.5em}

Since an even map is mapped to an even map from \eqref{1coc}, we consider
\[
    \begin{gathered}
        XB^2(L;M)_{\ev}:=\mathfrak{d}^1(XC^1(L;M)_{\ev})\subseteq XZ^2(L;M)_{\ev}.
    \end{gathered}
\]
\vspace{0.5mm}

As a result, we define
\begin{equation}\label{evco} H^2(L;M)_{\bar 0}=XZ^2(L;M)_{\bar 0}/XB^2(L;M)_{\bar 0}.\end{equation}

\vspace{0.5 em}

\subsection{Extensions of Lie superalgebras in characteristic 2.}\label{ext} Let $L$ be a Lie superalgebra in characteristic 2  and $M$ be an  abelian Lie superalgebra in  characteristic 2, that is $[M,M]=0$ and $s|_{M_{\od}}=0$.  An \emph{extension} of $L$ by $M$ is a short exact sequence 
\begin{equation}\label{eq:extension}
0 \longrightarrow M \xrightarrow{i} E \xrightarrow{\pi} L \longrightarrow 0,
\end{equation}
where \black $i$ and $\pi$ are \emph{morphisms} of Lie superalgebras.\\

\noindent Two such extensions $0\to M\xrightarrow{i_k}E_k\xrightarrow{\pi_k}L\to0,~ \text{for } k=1,2$
are called \emph{equivalent} if there exists an isomorphism $\sigma:E_1\to E_2$ such that the following diagram commutes:
\begin{equation}\label{eq:diagram}
  \xymatrix{
     && E_1 \ar[dr]^{\pi_1} \ar[dd]_{\sigma} & \\
     0 \ar[r] & M \ar[ur]^{i_1} \ar[dr]_{i_2} && L \ar[r] & 0. \\
     && E_2 \ar[ur]_{\pi_2} 
  }
\end{equation}

{\bf Technicalities.} Consider an  extension as given by \eqref{eq:extension}. Then, $i(M)$ is a homogeneous ideal of $E$, and $M$ has an $L$-module  structure. Indeed,  we have 
     \begin{align*}
    [\,i(m),x\,]&\in \mathrm{Ker}(\pi)=i(M), \quad \forall x \in L,\,\forall m\in M, \,\text{and}\\[4pt]
    s_E\circ i \, (m)&=i \circ s_M(m)\in i(M), \quad \forall m \in M_{\bar 1}.
    \end{align*}
  To check whether $i(M)$ is homogeneous, consider $c=c_{\ev}+c_{\od} \in i(M)$ where $c_{\ev} \in E_{\ev}$ and $c_{\od} \in E_{\od}$. We have $0 = \pi(c_0 +c_1) = \pi(c_0) + \pi(c_1)$. Since $\pi$ is even, we get  $\pi(c_0) =\black\pi(c_1)=0 \text{ which implies } c_0,\,c_1\in\mathrm{Ker}(\pi)=i(M). \text{ Hence, } c_0\in i(M)_{\bar 0} \text{ and } c_1 \in i(M)_{\bar 1}$, and thus, $i(M)$ is homogeneous.

  Moreover, we define an $L$-module structure on $M$ by 
  \begin{equation}\label{eq:module}
    x\cdot m = i^{-1}\bigl([\,\tilde{x},i(m)\bigr]\bigr) \in i(M),\quad \forall x\in L,\ \forall m\in M.
  \end{equation}
  where $\pi(\tilde{x})=x, \, \text{ for some } \,\tilde{x} \in E$ by the surjectivity of $\pi$. In the case where $i(M)\subseteq \mathfrak{z}(E)$, this action is identically zero and $M$ becomes a trivial $L$-module. These extensions are called \emph{central extensions}.
  \vspace{1mm}

 We now show that $H^2(L;M)_{\ev}$ captures all central extensions of Lie  superalgebras by abelian Lie superalgebras. 

\sssbegin{Theorem} \label{thm1}
Let $L$ be a Lie superalgebra and $M$ be an abelian Lie superalgebra. Then, the equivalence classes of central extensions of $L$ by $M$ are are classified by $H^2(L;M)_{\bar 0}$.
\end{Theorem}

\begin{proof}
The proof follows five steps.

\vspace{1 mm}

\noindent\underline{Step 1.} Consider a central extension of $L$ by $M$ as in \eqref{eq:extension}. Then, let $\rho:L \longrightarrow E$ be an even, linear section of the surjective morphism $\pi$, such that $\pi \circ \rho=id_{L}$. One can deduce that 
\begin{equation*}
[\rho(x),\rho(y)]-\rho([x,y]) \in \mathrm{Ker}(\pi) \text{ and } s\circ\rho(z)-\rho\circ s(z) \in \mathrm{Ker}(\pi), \ \ \forall \, x,y\in L, z \in L_{\od}.
\end{equation*}

\vspace{0.5 em}

\noindent We define the following two maps:
\begin{align}\label{11}
        \varphi: L\times L&\longrightarrow M,~ (x,y)\longmapsto i^{-1}\Bigl([\rho(x),\rho(y)]-\rho([x,y])\Bigl);\\
        \omega: L_\ev&\longrightarrow M_\ev,~x\longmapsto i^{-1}\Bigl(s\circ\rho(x)-\rho\circ s(x)\Bigl).
\end{align}

\noindent We now show that $(\varphi,\omega)\in XZ^2(L;M)_{\ev}.$

\vspace{0.5 em}

\noindent The fact that the $2$-cocycle is even follows from \eqref{cograd} and \eqref{11}. Moreover, $\varphi$ is clearly bilinear and skew symmetric, and $\omega$ a squaring. We show $\varphi(x,y)=\omega(x+y)+\omega(x)+\omega(y), \ \forall x,y \in L_{\od}.$ Indeed,
    \begin{align*}
        \omega(x+y)&=i^{-1}\bigl(s\circ\rho(x+y)-\rho\circ s(x+y)\bigr)\\
        &=i^{-1}\bigl(s(\rho(x)+\rho(y))-\rho([x,y]+s(x)+s(y)\bigr)\\
        &=i^{-1}\bigl([\rho(x),\rho(y)]+s\circ\rho(x)+s\circ\rho(y) -\rho([x,y])-\rho\circ s(x)- \rho \circ s(y)\bigr)\\
        &= \omega(x)+\omega(y)+\varphi(x,y). 
    \end{align*}
    
\vspace{0.5 em}

\noindent We verify that $\mathfrak{d}^2(\varphi,\omega) = 0$ (see Definition \eqref{2coc}). We show first that $d_{\ce}^2(\varphi)(x,y,z)=0, \ \forall x,y,z \in L$. Since the action is trivial we have 
\[d_{\ce}^2(\varphi)(x,y,z)= \underset{x,y,z}{\circlearrowleft}i^{-1}\bigl([\rho([x,y]),\rho(z)]-\rho([[x,y],z)]\bigr)\]
Since $[\rho(x),\rho(y)]-\rho([x,y])\in i(M)$ and $[i(M),E]=0$ (as $i(M)\subseteq \fz(E)$), we have $[[\rho(x),\rho(y)]-\rho([x,y]),\rho(z)]=0$. 

Thus, we have $[[\rho(x),\rho(y)],\rho(z)]=[\rho([x,y]),\rho(z)]$. This, together with the linearity of $\rho$ and $i^{-1}$ gives us:
\[d_{\ce}^2(\varphi)(x,y,z)= i^{-1}\bigl(\, \underset{x,y,z}{\circlearrowleft}([[\rho(x),\rho(y)],\rho(z)])-\rho(\,\underset{x,y,z}{\circlearrowleft}[[x,y],z])\bigr).\]
By the Jacobi Identity, we get $d^2(\varphi)(x,y,z)=0$. 

\vspace{0.5 em}

We now show that $\delta^2(\varphi,\omega)(x,z)=0,\, \forall\,x\in L_{\od},\, \forall\,z\in L$. Since the action is trivial we have
\vspace{1mm}
\begin{equation}\label{13}
\delta^2(\varphi,\omega)(x,z)=i^{-1}\bigl([\rho \circ s(x),\rho(z)]-\rho([s(x),z])+[\rho([x,z]),\rho(x)]-\rho([[x,z],x]\bigl).
\vspace{1mm}
\end{equation}

 Since $\rho \circ s -s\circ\rho\in i(M)$ and $[i(M),E]=0$, we have \([\rho \circ s(x), \rho(z)]=[s \circ \rho(x),\rho(z)]\). Hence, we rewrite \eqref{13} as
\begin{equation*}
\delta^2(\varphi,\omega)(x,z)=i^{-1}\bigl([s \circ \rho(x),\rho(z)]-\rho([s(x),z])+[\rho([x,z]),\rho(x)]-\rho([[x,z],x]\bigr).
\end{equation*}

Using $[s(x),z]=[x,[x,z]]$, this reduces to
    \begin{align*}
        \delta^2(\varphi,\omega)(x,z)&=i^{-1}\bigl([s \circ \rho(x),\rho(z)]+[\rho([x,z]),\rho(x)]\bigr)\\
        &=i^{-1}\bigl([\rho(x),[\rho(x),\rho(z)]]+[\rho([x,z]),\rho(x)]\bigr)    
    \end{align*}
Similarly, since $[\rho(x),\rho(z)]-\rho([x,z])\in i(M)$ and $[i(M),E]=0$, we have $[[\rho(x),\rho(z)]-\rho([x,z]),\rho(x)]=0$. Hence, $\delta^2(\varphi,\omega)=0$. Therefore, $(\varphi,\omega)\in XZ^2(L;M)_{\ev}$.

\vspace{1 em}

\noindent\underline{Step 2.} Let $(\varphi,\omega) \in XZ^2(L;M)_{\ev}.$ We define the following bracket and squaring maps, respectively, on $ L \oplus M$:
\begin{align}\label{14}
    [x+m,y+n]_{\mathrm{New}} &:= [x,y]_L+\varphi(x,y), \ \ \forall x,y\in L,\ \forall m,n \in M;\\
    s_{ \mathrm{New}}(x+m) &:= s_{L}(x) + \omega(x), \ \ \forall x\in L_{\od},\ \forall m \in M_{\od}.
\end{align}

We show that this defines a valid Lie superalgebra in characteristic 2 on $L\oplus M$. To check whether \eqref{14} imposes a valid Lie algebra structure on $(L \oplus M)_{\ev}$ is routine. We check if the odd part is a symmetric module of the even part. Let $x+m \in (L \oplus M)_{\ev}$ and $y+n \in (L \oplus M)_{\od}$. Then,
    \begin{align*}
        [x+m, y+n]_{N} &= [x,y]_{L}+ \varphi(x,y)=[y,x]_{L}+\varphi(y,x)=[y+n, x+m]_{N}.
    \end{align*}

We now check the properties of the squaring $s_{\mathrm{New}}$. Let $x+m, y+n \in (L \oplus M)_{\od}$. Then, we have 
        \begin{align*}
            [x+m,y+n]_{N}&=[x,y]_{L}+\varphi(x,y)\\
            &=s_L(x+y)+s_L(x)+s_L(y)+\omega(x+y)+\omega(x)+\omega(y)\\
            &=s_N((x+m)+(y+n))+s_N(x+m)+s_N(y+n).
        \end{align*}
Now, let $x+m \in (L\oplus M)_{\od}$ \,, $y+n \in L \oplus M$. Then, we have 
        \begin{align*}
            [x+m,[x+m,y+n]_{N}]_{N}&=[x+m,[x,y]_{L}+\varphi(x,y)]_{N}\\
            &=[x,[x,y]_{L}]_{L}+\varphi(x,[x,y]_{L})
        \end{align*}
        And,
        \begin{align*}
            [s_{N}(x+m), y+n]_{N}&=[s_{L}(x)+\omega(x), y+n]_{N}\\
            &=[s_L(x),y]_{L}+\varphi(s_L(x),y)\\
            &=[x,[x,y]_L]_L+\varphi(x,[x,y]_L) \qquad \mathrm{(from} \ \delta^2(\varphi,\omega)=0 \text{ in  \eqref{2coc})}\\
        \end{align*}
Thus, $[s_{N}(x+m), y+n]_{N}=[x+m,[x+m,y+n]_{N}]_{N}$.

\vspace{1em}

\noindent\underline{Step 3.} Let $E_{1}$ and $E_{2}$ be two equivalent central extensions of $L$ by $M$. We denote by $\sigma:E_1 \longrightarrow E_2$, the isomorphism of Lie superalgebras that realizes the equivalence in the initial commutative diagram \eqref{eq:diagram}. We show that $(\varphi_1,\omega_1)$ and $(\varphi_2,\omega_2)$ (which are the $2$-cocycles derived from $E_1$ and $E_2$ respectively using Step 1) are cohomologous, that is $(\varphi_2-\varphi_1, \omega_2-\omega_1) \in XB^2(L;M)_{\ev}$. Denote by $\rho_1$ (respectively $\rho_2$) the section of $\pi_1$ (respectively $\pi_2$). We have $\pi_2(\sigma\rho_1-\rho_2)=0$, and therefore we can define an even map $\psi:L\longrightarrow M, \ x\longmapsto i_2^{-1}\circ(\sigma\rho_1-\rho_2)(x)$. Let us now consider the following map $i_2 \circ\psi:L\longrightarrow i_2(M) \cong M, \ x\longmapsto(\sigma\rho_1-\rho_2)(x)$. We attempt to show first that $\varphi_2(x,y)=\varphi_1(x,y)+d_{\ce}^1(\psi(x,y)), \ \forall x,y \in L$.

\vspace{0.5em}

Recall that, $[\rho_k(x),\rho_k(y)]-\rho_k([x,y])=i_k \circ\varphi_k(x,y)$ for $k=1,2$. Then, we have 
\begin{align*}
\sigma([\rho_1(x),\rho_1(y)])=[\sigma\rho_1(x),\sigma\rho_1(y)]=[\rho_2(x)+i_2\psi(x),\rho_2(y)+i_2\psi(y)].
\end{align*}

Since $\mathrm{Im}(i_2\circ\psi) \subseteq i_2(M)\subseteq \mathfrak{z}(E)$, we get  \[
[i_2\psi(x),\rho_2(y)]=[\rho_1(x),i_2\psi(y)]=[i_2\psi(x),i_2\psi(y)]=0.
\]
This gives us $\sigma[\rho_1(x),\rho_1(y)]=[\rho_2(x),\rho_2(y)].$ Moreover, we also have
\begin{align*}
    \sigma([\rho_1(x),\rho_1(y)]) &= \sigma([\rho_1([x,y])+i_1\varphi_1(x,y)) \black\\
    &=\sigma\rho_1([x,y]) + \sigma(i_1\varphi_1(x,y))\\
    &=\rho_2([x,y]) +i_2\psi([x,y]) +i_2\varphi_1(x,y) 
\end{align*}

Together with $[\rho_2(x),\rho_2(y)]=\rho_2([x,y])+i_2\varphi_2(x,y)$ and the injectivity of $i_2$, we get \(\varphi_2=\varphi_1+d^1_{\ce}\psi.\)

\vspace{0.5em}
Next, we verify {$\delta^1(\psi)(x)= \omega_2(x)-\omega_1(x), \ \forall x \in L_{\od}$}. Then, similar to the previous verification, we have 
\begin{align*}
    \sigma(s \circ \rho_1(x))&=\sigma(\rho_1 \circ s(x) + i_1\omega_1(x))\\
    &=\sigma(\rho_1 \circ s(x)) +i_2\omega_1(x))\\
    &=\rho_2\circ s(x)+i_2\psi(s(x)) + i_2\omega_1(x).
\end{align*}
But
\begin{align*}
\sigma(s\circ \rho_1)(x)&=s(\sigma \circ \rho_1)(x)\\
&=s(\rho_2(x)+i_2\psi(x))\\
&=[\rho_2(x),i_2\psi(x)]+s\circ \rho_2(x) +s(i_2\psi(x))\\
&=s\circ \rho_2(x) \qquad (\mathrm{Recall} \ s(M_{\od})=0)\\
&=\rho_2 \circ s(x)+i_2\omega_2(x).
\end{align*}

Equating the two gives us $\omega_2 = \delta^1\psi +\omega_1$. For sums in particular we elucidate the following, though it already follows from the previous verification. Let $x,y \in L_{\od}$ where $x$ and $y$ are basis elements. Then,
\begin{align*}
(\omega_2-\omega_1)(x+y)&=(\varphi_2-\varphi_1)(x,y)+(\omega_2-\omega_1)(x)+(\omega_2-\omega_1)(y)\\
&=d^1_{\ce}\psi(x,y)+\delta^1\psi(x)+\delta^1\psi(y)\\
&=\psi([x,y]+s(x)+s(y))\\
&=\psi(s(x+y))\\
&=\delta^1\psi(x+y).
\end{align*}

Hence, $(\varphi_2-\varphi_1, \omega_2-\omega_1) \in XB^2(L;M)_{\ev}.$

\vspace{1 em}

\noindent\underline{Step 4.} Let $\psi: L\longrightarrow M$ be an even linear map and $(\varphi_1,\omega_1)\in XZ^2(L;M)_{\ev}$. Define 
\[
(\varphi_2,\omega_2)=(\varphi_1,\omega_1)+\mathfrak{d}^1(\psi).
\]
The final step of this proof aims to show that the central extensions $E_1$ and $E_2$, given by $(\varphi_1,\omega_1)$ and $(\varphi_2,\omega_2)$ respectively are equivalent. Consider the bijective linear map $\sigma\colon E_1\longrightarrow E_2$ defined by 
\[
\sigma(x+m)=x+m+\psi(x).
\]
This map makes the diagram \eqref{eq:diagram} commute. It remains to show that $\sigma$ is a morphism of Lie superalgebras. Let $x_1+m_1$ and $x_2+m_2$ be elements of $E_1$.  Then, the commutativity with the bracket follows from the proof in Step 4 in \cite[Theorem 3.5.3]{BE}. 

\vspace{0.5 em}

\noindent Now, let $x+m$ be an odd element of $E_1$. Then, we have
\begin{align*}
    s(\sigma(x+m))&=s(x+m+\psi(x))\\
    &=s_L(x)+\omega_2(x)\\
    &=s_L(x)+\omega_1(x)+\psi(s_L(x))\\
    &=\sigma(s_L(x)+\omega_1(x))\\
    &=\sigma(s(x+m)).
\end{align*}
Thus, $E_1$ and $E_2$ are equivalent.

\sssbegin{Remark}
We note here that up until now we have defined two maps in Step 1 and Step 2 that assign a cocycle to each central extension and vice versa, respectively. Steps 3 and 4 show that these maps are well-defined over equivalence classes. To prove bijectivity, however, we need to prove that both these maps are 2-sided inverses of one another. This is done formally in Step 5, and is missing from \cite[Theorem 3.5.3]{BE}. However, due to the nature of the proof for the following step, one verifies that it can indeed be adapted to the restricted case, and hence, \cite[Theorem 3.5.3]{BE} holds.
\end{Remark}

\vspace{0.5 em}

\noindent\underline{Step 5.} Since we have proven that the methods in Step 1 and Step 2 are well-defined over equivalence classes, we now prove that both these processes are inverses of one another. Formally, if we let 
\begin{align*}
\Phi \colon &\mathcal{E} \longrightarrow H^2(L;M)_{\ev}\\
&[E] \longmapsto [(\varphi,\omega)]
\end{align*}
be the map that assigns a cohomology class to each equivalence class of central extensions by Steps 1 and 3, and we let
\begin{align*}
\Psi \colon &H^2(L;M)_{\ev} \longrightarrow \mathcal{E}\\
&[(\varphi,\omega)] \longmapsto [E]
\end{align*}
be the map that assigns an equivalence class of central extensions to each cohomology class by Steps 2 and 4, then we prove that $\Phi$ and $\Psi$ are two-sided inverses of one another. 

\vspace{0.5 em}

Let us first prove that $\Phi \circ \Psi = \id_{H^2(L;M)_{\ev}}$. Let $(\varphi, \omega) \in H^2(L;M)_{\ev}$.  We build the central extension $E=L \oplus M$ using Step 2. Now, we pick an even, linear section $\rho$. Since $\pi$ is just the projection map from our construction in Step 2 and $\pi \circ \rho = \mathrm{id}_{L \oplus M}$, we have that $\rho(x)=x+\psi(x)$ for an arbitrary linear map $\psi \colon L \longrightarrow M$.

We define the cocycle $(\phi, \gamma)$ as in Step 1. We have 
$\phi(x,y)=i^{-1}\Bigl([\rho(x),\rho(y)]_{E} -\rho([x,y]_{L})\Bigl)$ and $\gamma(z)=i^{-1}\Bigl(s_{E}\circ\rho(z)-\rho\circ s_{L}(z)\Bigl).$ Since $i$ is just the inclusion map $i(m)=0+m$ from our construction in Step 2, we retain the same equations while dropping $i^{-1}$. Then, we have
\begin{align*}
\phi(x,y)&=[\rho(x),\rho(y)]_{E} -\rho([x,y]_{L})\\
&=[x+\psi(x), y+\psi(y)]_{E}-\rho([x,y]_{L})\\
&=[x,y]_L+\varphi(x,y)-[x,y]_{L}-\psi([x,y])\\
&=\varphi(x,y)-\psi([x,y])\\
&=\varphi(x,y)+d_{\ce}^1\psi(x,y).
\end{align*}

Similarly,
\begin{align*}
\gamma(z)&=s_{E}\circ \rho(z)-\rho \circ s_{L}(z)\\
&= s_{E}(z+\psi(z))-\rho \circ s_{L}(z)\\
&=s_L(z)+\omega(z)-s_L(z)-\psi(s_L(z))\\
&=\omega(z)+\delta^1\psi(z).
\end{align*}

Hence, we see that $(\phi,\gamma)$ is cohomologous to the original cocycle used to build the extension $(\varphi, \omega)$, and we have proven invertibility from one side.

\vspace{0.5 em}

Now, we prove $\Psi \circ \Phi = \id_{\mathcal{E}}$. Let $E$ be an arbitrary central extension of the form as in \eqref{eq:extension}. Then, by Step 1 we produce a $2$-cocycle.  Now, using this 2-cocycle we construct a central extension $E'=L\oplus M$ as in Step 2. Explicitly, we have
\begin{equation*}
0 \longrightarrow M \xrightarrow{i_N} L \oplus M \xrightarrow{\pi_N} L \longrightarrow 0,
\end{equation*}
where $i_N(m)=m$ and $\pi_{N}(x+m)=x$ for all $x \in L, m\in M$ strictly. We continue to denote by $i$ and $\pi$ the morphisms from the original extension $E$. We now try to prove an equivalence of these extensions. Let $\rho: L \longrightarrow E$ be an even, linear section of $\pi$ such that $\pi \circ \rho=\id_{L}$. Any $z\in E$ can be written as $z= \rho\circ\pi(z)+(z-\rho\circ\pi(z))$. If we set $x=\pi(z) \in L$ and see that $z-\rho(x) \in i(M)$ (as $\pi(z-\rho(x)) \in \Ker(\pi))$, we have the  decomposition $z=\rho(x)+i(m)$. Hence, we define the bijective linear map $\sigma: L \oplus M \longrightarrow E$ given by $\sigma(x+m)=\rho(x)+i(m)$. Indeed, one deduces that the map $\sigma^{-1}(z)=\pi(z)+i^{-1}(z-\rho\circ\pi(z))$ is a valid 2-sided inverse. We would now like to prove that the following commutative diagram holds
\begin{equation}\label{diagram2}
  \xymatrix{
     && L \oplus M \ar[dr]^{\pi_N} \ar[dd]_{\sigma} & \\
     0 \ar[r] & M \ar[ur]^{i_N} \ar[dr]_{i} && L \ar[r] & 0, \\
     && E \ar[ur]_{\pi} 
  }
\end{equation}
where $\sigma$ is an isomorphism.

We first prove that $\sigma$ defines an isomorphism. We first check for the Lie bracket. Let $x_1+m_1, x_2+m_2 \in L \oplus M$. Then,
\begin{align*}
\sigma([x_1+m_1,x_2+m_2]_{L\oplus M})&=\sigma([x_1,x_2]_L+\varphi(x_1,x_2))\\
&=\rho([x_1,x_2]_L)+i(\varphi(x_1,x_2))\\
&=\rho([x_1,x_2]_L)+[\rho(x_1),\rho(x_2)]_E-\rho([x_1,x_2]_L)\\
&=[\rho(x_1),\rho(x_2)]_E.
\end{align*}

On the other hand,
\begin{align*}
[\sigma(x_1+m_1),\sigma(x_2+m_2)]_E&=[\rho(x_1)+i(m_1),\rho(x_2)+i(m_2)]_E\\
&=[\rho(x_1),\rho(x_2)]_E.
\end{align*}
We now move on to the squaring. Let $x+m \in (L \oplus M)_{\od}$. Then,
\begin{align*}
\sigma(s(x+m)_{L \oplus M})&=\sigma(s_L(x)+\omega(x))\\
&=\rho(s_L(x))+i(\omega(x))\\
&=\rho(s_L(x))+s_E(\rho(x))-\rho(s_L(x))\\
&=s_E(\rho(x)).
\end{align*}

Similarly, we also have
\begin{align*}
s(\sigma(x+m))_{E}&=s(\rho(x)+i(m))_E\\
&=[\rho(x),i(m)]_E+s(\rho(x))_E+s(i(m))_E\\
&=s(\rho(x))_E.
\end{align*}

Hence, $\sigma$ is an isomorphism. Finally, we see that $\sigma$ makes the diagram \eqref{diagram2} commute rather straightforwardly. Indeed, we have
\[\sigma\circ i_N(m)=\sigma(m)=i(m).\] And,
\[\pi\circ \sigma(x+m)=\pi(\rho(x)+i(m))=x=\pi_N(x+m).\] 

Hence, we see that $E'=L \oplus M$ is equivalent to the original extension $E$, and we have proved the invertibility from the other side.
\end{proof}

\sssbegin{Corollary}\label{cor}Let $L$ be a Lie superalgebra. Then, the equivalence classes of one-dimensional central extensions of $L$ by $M$ are classified by $H^2(L;\mathbb{K})_{\ev}.$ 
\end{Corollary}

\section{Inductive Theorem for Nilpotent Lie superalgebras}\
This section aims to show that any $n$-dimensional nilpotent Lie superalgebra can be obtained as a central extension of an $(n-1)$-dimensional nilpotent Lie superalgebra using Corollary \ref{cor}. The method we use is adapted and built on the approach taken in \cite[Section 5.1]{BE} which requires a preliminary set-up.

\subsection{Equivalent Scalar Cocycles.}\label{ModStr} Let $L$ be a Lie superalgebra, and let $\Aut(L)$ denote the the automorphism group of $L$. Let $M$ be an $L$-module. We impose an $\Aut(L)$-module structure on $XZ^2(L;M)$ as follows. Let $(\varphi,\omega)\in XZ^2(L;M)$ and $A\in \Aut(L)$. Then, we define the action of $\Aut(L)$ on $XZ^2(L;M)$, given by $A\cdot(\varphi,\omega):=(A\varphi,A\omega)$, as

\begin{equation}\label{action}
    \begin{cases}
        (A\varphi)(x,y) &= \varphi(A(x),A(y)), \quad  \forall x,y \in L; \\
(A\omega)(x) & = \omega(A(x)), \quad  \forall x \in L_{\od}. 
    \end{cases}
\end{equation}

It is straightforward to see that the coboundary group $XB^2(L;M)$ is an $\Aut(L)$-submodule of $XZ^2(L;M)$, and hence we may define the $\Aut(L)$-module structure on $H^2(L;M)$. We say two cocycles are \emph{equivalent} if they belong to the same $\Aut(L)$-orbit. 

\vspace{0.5 em}

 We now show that any two equivalent cocycles belonging to $H^2(L;\K)$ produce isomorphic central extensions. Unlike in \cite[Proposition 5.2.1]{BE} where, owing to Jacobson's Theorem \cite[Theorem 2.2.1]{BE}, the proof could be reduced to only considering the Lie bracket, we must also verify the compatibility of any potential isomorphism with the squaring. 

\sssbegin{Proposition}\label{prop1}Let $(L,[\cdot,\cdot],s)$ be a nilpotent Lie superalgebra and let $(\varphi,\omega)$ and $(\psi,\gamma)$ be two equivalent scalar $2$-cocycles. Let $(L_{\varphi},[\cdot,\cdot]_{\varphi},s_{\varphi})$ and $(L_{\psi},[\cdot,\cdot]_{\psi},s_{\psi})$ be two central extensions of $L$ by an element $X$ such that $|X|=|\varphi|=|\psi|$, with brackets given by
\[
[\cdot,\cdot]_{\varphi} := [\cdot,\cdot] + \varphi(\cdot,\cdot)X, \quad [\cdot,\cdot]_{\psi} := [\cdot,\cdot] + \psi(\cdot,\cdot)X,
\]
and the squarings given by
\[
s_{\varphi}(\cdot):=s(\cdot) + \omega(\cdot)X, \quad s_{\psi}(\cdot):=s(\cdot) + \gamma(\cdot)X
\]
respectively. Then, the Lie superalgebras $L_{\varphi}$ and $L_{\psi}$ are isomorphic.
\end{Proposition} 

\begin{proof} Since $\varphi$ and $\psi$ are equivalent, this implies that there must exist an $A \in \Aut(L)$ such that
\begin{align*}
\varphi(x,y) &= \psi(Ax,Ay) \quad \forall x,y \in L,\\
\omega(x)&=\gamma(Ax) \quad \forall x,y\in L_{\od}.
\end{align*}

\noindent We define the map $f : L_{\varphi} \to L_{\psi}$ with $f(x) = Ax, \ \forall x \in L \text{ such that } f(X) = X$. It is straightforward to verify that $f$ is both bijective and linear. Let $x,y \in L$. Then, the commutativity with the bracket follows from \cite[Proposition 5.2.1]{BE}.

\vspace{0.5 em}

\noindent For $x\in L_{\od}$, we have 
\[
\begin{aligned}
f\circ s_\varphi(x) &= f(s(x)+\omega(x)X) \\
&= A\circ s(x) + \omega(x)f(X) \\
&= s(A(x)) + \gamma(Ax) \\
&=s_{\psi}\circ f(x).
\end{aligned}
\]

\noindent Hence, $f$ defines an isomorphism, and $L_{\varphi}$ and $L_{\psi}$ are isomorphic Lie superalgebras.
\end{proof}

\vspace{0.5 em}

\subsection{Inductive Theorem.} We now move on to the main theorem we shall use for the classification. 

\sssbegin{Theorem}\label{indthm}
Let L be a nilpotent Lie superalgebra of dimension $n$. Then, $L$ is isomorphic to a central extension by a scalar 2-cocycle of a nilpotent Lie superalgebra of dimension $n-1$.
\end{Theorem}
\begin{proof}
By assumption, nilpotence implies that $\fz(L)\neq 0$. We then split the proof into two cases, as follows.

\vspace{0.5 em}

\noindent\underline{$\mathfrak{z}(L)_{\ev} \neq 0$.} Let $Y\in\mathfrak{z}(L)_{\ev}$ and let $I$ be the one-dimensional, exclusively even ideal generated by $Y$. Then, $I$ is abelian as $Y$ belongs to the center. Then, we define the quotient Lie superalgebra $H=L/I$. We consider the following short exact sequence of Lie superalgebras 
\[0\longrightarrow I\longrightarrow L\overset{\pi}{\longrightarrow} H\longrightarrow0, \]
where $\pi: L \longrightarrow H$ is the projection map, which is a Lie superalgebra morphism. We repeat the construction from Theorem \ref{thm1} to produce the 2-cocycle
\begin{align*}
    \varphi:H\times H \longrightarrow I,\ & (x,y)\longmapsto [\rho(x),\rho(y)]-\rho([x,y]),\\
    \omega: H_{\od}\longrightarrow I_{\ev},\ & x\longmapsto s\circ\rho(x)-\rho\circ s(x),
\end{align*}

\noindent where $(\varphi,\omega)\in XZ^{2}(H;I)_{\ev}.$ 

\vspace{0.5 em}

Since $\mathrm{Im}(\varphi), \mathrm{Im}(\omega)\subseteq\mathrm{Ker}(\pi)=I$, there exists a pair $(\Delta,\gamma)$ where $\Delta:H \times H\longrightarrow \mathbb{K}$ is an even bilinear form such that $\varphi(x,y)=\Delta(x,y)Y, \ \forall x,y\in H$, and $\gamma:H_{\od}\longrightarrow \mathbb{K}$ is an odd scalar map such that $\omega(x)=\gamma(x)Y, \ \forall x \in H_{\od}$. We already have $(\varphi,\omega)\in XZ^2(H;I)_{\ev}$ from Step 2 of Theorem \ref{thm1}. It follows that $(\Delta,\gamma)\in XZ^2(H;\mathbb{K})_{\ev}$. We have that $L$ is isomorphic to the one-dimensional central extension of $H$ by the even scalar 2-cocycle $(\Delta,\gamma)$ by Step 5 of Theorem \ref{thm1}.

\noindent\underline{$\mathfrak{z}(L)_{\ev}=0$.} In this case, the center is purely odd. Let $Y\in\mathfrak{z}(L)_{\od}$ and let $I$ be the one-dimensional purely odd ideal generated by $Y$. Clearly, $I$ is abelian as $Y$ belongs to the center. Then, we define the quotient Lie superalgebra $H=L/I$. We again have the same short exact sequence of Lie superalgebras as in the first case, and we repeat the construction of $(\varphi,\omega)\in XZ^2(H;I)_{\ev}$. However, in this case, since $I$ is purely odd, $\omega\equiv 0$.

Similar to the previous step, we derive the pair $(\Delta,0) \in Z^2(H;\mathbb{K})_{\od}$, and prove that $L$ is isomorphic to the central extension produced by $(\Delta,0)$. 
\end{proof}

\section{Classification in Dimension 3} The classification method works as follows. 

\subsection{Equivalence Classes of Scalar 2-cocycles.}\label{method}We begin by first determining the equivalence classes of the the scalar $2$-cocycles for every Lie superalgebra. Again, we recall that two scalar $2$-cocycles $(\varphi,\omega),(\psi,\gamma)\in XZ^2(L;\K)$ of a Lie superalgebra $L$ are \emph{equivalent} if there exists an automorphism of $L$, denoted $A\in \text{Aut}(L)$, such that $(\varphi,\omega) = A \cdot (\psi,\gamma)$, given by
\begin{equation}
\begin{cases}\label{conditions}
\varphi(x,y) &= \psi(Ax,Ay) \quad \forall x,y \in L,\\
\omega(x) & =\gamma(Ax) \quad \forall x\in L_{\od}.
\end{cases}
\end{equation}

Let $L$ be  a nilpotent Lie superalgebra of superdimension $(n|m)$. Then, we denote the basis elements of $L$ by $\{e_1,...,e_{n+m}\}$, where $L_{\ev}= \text{Span}\{e_1,...,e_n\}$ and $L_{\od}= \text{Span}\{e_{n+1},...,e_{n+m}\}$.

The basis elements of the scalar $2$ cocycles are given by the elements
\begin{equation*}
\begin{cases}
\Delta_{i,j}:L\times L\longrightarrow \K, &~1\leq i \leq   j\leq n+m;\\
\omega_{k}: L_{\od}\longrightarrow\K, &~n+1\leq k \leq n+m. 
\end{cases}
\end{equation*}
where $\Delta_{i,j}(e_s,e_t)=\delta_{i,s}\delta_{t,j}$ and $\omega_{k}(e_l)=\delta_{k,l}$. We crucially note that $\Delta_{i,j}=\Delta_{j,i}$, and that $\Delta_{i,j}\equiv0$ for all $i=j$ from the definition of a $2$-cocycle $\in XZ^2$ in \ref{cohom}. \footnote{We would also like to reiterate here that the squaring is not linear, and a $2$-cocycle of the form $(\varphi,\omega)$ must satisfy $\omega (\lambda x)= \lambda^2 \omega  (x)$ and $\omega(x+y)=\omega(x)+\omega(y)+\varphi(x,y)$.}

For our convenience, we reduce the notation above and let $\Delta_{i,j}$ denote the entire class of maps  $\Delta_{i,j}(e_i,e_j)=\lambda_{ij} \in \K\setminus\{0\}$, and $0$ otherwise. Similarly $\omega_{k}(e_k)=\mu_{k}\in\K\setminus\{0\}$, and $0$ otherwise. An important point to mention here is that if the first part of the two cocycle is $0$, then indeed, any input is mapped to $0$. This is not the case for the squaring, where the $0$-squaring may send the sum of the basis elements to a non-zero scalar value if the first part of 2-cocycle is non-zero on the basis elements (see Section \ref{cohom}). This is crucial while verifying the conditions in \eqref{conditions}.

To compute the $2$-cocycles themselves, we use the Mathematica package \texttt{SuperLie} (see \cite{Kr}). We then determine the equivalence classes by computing the automorphisms of $L$ that make the scalar $2$-cocycles equivalent. For every Lie superalgebra, we begin by first stating the necessary conditions for each automorphism to in fact be a valid Lie superalgebra morphism. Then, we derive the conditions that the relations in \eqref{conditions} impose on the basis elements and aim to satisfy them while maintaining $\text{det}(A)\neq0$ (the latter being a condition we always check even if it is not explicitly mentioned). If $(\varphi, \omega)$ and $(\varphi', \omega')$ are equivalent, we denote it as $(\varphi, \omega)\cong(\varphi', \omega')$. We determine that two $2$-cocyles can not be equivalent when we reach an unrectifiable contradiction that prevents the relations in \eqref{conditions}.

\sssbegin{Remark}
With regards to notation, we differentiate the values of the 2-cocycles we wish to prove the automorphism for by writing the latter one with a $``~'~"$. For example, when writing $(\Delta_{1,2},0)$ is equivalent to $(\Delta_{1,3},0)$, we write $(\Delta_{1,2},0)\cong(\Delta_{1,3}',0)$, where $\Delta_{1,3}'(e_1,e_3)=\lambda_{13}' \in \K \setminus \{0\}$. 
\end{Remark}

\subsection{Automorphism Computation for dim(\texorpdfstring{$L$}{L})=2.} We  compute the scalar $2$-cocycles for each of the $2$-dimensional nilpotent Lie superalgebras in \ref{nilclassif2} using \texttt{SuperLie} (see \cite{Kr}). 

\sssbegin{Lemma}
Suppose $L$ is a nilpotent Lie superalgebra with $\dim(L)=2$ over an algebraically closed field of characteristic 2. Then, the scalar $2$-cocycle basis of $L$ is given by
\end{Lemma}

{\renewcommand{\arraystretch}{1.5}
\small
\begin{center}
    \makebox[\textwidth][c]{
    \begin{tabular}{|c|c|c|c|}
    \hline
          \textbf{Name} &\textbf{sdim} & \textbf{Even Cocycles}& \textbf{Odd Cocycles} \\\hline
          $\mathbf{L^{1}_{0|2}}$  & $(0|2)$ & $(\Delta_{1,2}), ~(0, \omega_1), ~(0, \omega_2)$ & --- \\\hline
          $\mathbf{L^{2}_{1|1}}$  & $(1|1)$ & --- & --- \\\hline
          $\mathbf{L^{3}_{1|1}}$  & $(1|1)$ & $(0,\omega_2)$ & $(\Delta_{1,2},0)$ \\\hline
          $\mathbf{L^{1}_{2|0}}$  & $(0|2)$ & $(\Delta_{1,2})$ & --- \\\hline
          
    \end{tabular}
    }
   \\~\\ Scalar $2$-cocyle basis of nilpotent Lie superalgebras of dimension $2$.
\end{center}
}
\normalsize

\vspace{0.2 em}

\sssbegin{Theorem}\label{4.2.1} Suppose $L$ is a nilpotent Lie superalgebra with $\dim(L)=2$ over an algebraically closed field of characteristic $2$. Then, the equivalence classes of scalar $2$-cocycles on L are given by

\begin{itemize}
\item[] $L=\mathbf{L^{1}_{0|2}}: (0,0) , ~(\Delta_{1,2},0), ~(0,\omega_1) $

\item[] $L=\mathbf{L^{2}_{1|1}}: (0,0) $

\item[] $L=\mathbf{L^{3}_{1|1}}: (0,0),~(\Delta_{1,2},0),~(0,\omega_2) $

\item[] $L=\mathbf{L^{2}_{2|0}}: (0,0),~(\Delta_{1,2},0) $

\end{itemize}

\end{Theorem}

\begin{proof}
For the first case $L=\mathbf{L^{1}_{0|2}}$, the general form of the automorphism $A\in \text{Aut}(L)$ is given by
\begin{equation*}
Ae_1=a_1e_1+a_2e_2, ~Ae_2=b_1e_1+b_2e_2.
\end{equation*}

We have dim($H^2(L;\K))=3$ with the basis elements being $\{(\Delta_{1,2},0),(0,\omega_1),(0,\omega_2)\}$, all of which are even. Indeed we find almost immediately that $(\Delta_{1,2},0)\ncong (0,\omega_1')$, finding a contradiction with the first condition in \eqref{conditions}. We find that $(0,\omega_2)\cong (0,\omega_1')$ and $(0,\omega_1)+(0,\omega_2) \cong (0,\omega_1')$ using the automorphisms $Ae_1=e_2, ~Ae_2=\sqrt{\frac{\mu_2}{\mu_1'}}e_1$ and $Ae_1=\sqrt{\frac{\mu_1}{\mu_1'}}e_1, ~Ae_2=\sqrt{\frac{\mu_2}{\mu_1'}}e_1+e_2$, respectively.

We also find that $(\Delta_{1,2},0)+(0,\omega_1)\cong ({\Delta}_{1,2}',0)$ with $Ae_1=\frac{\mu_1}{\lambda_{12}'}e_1+e_2, ~Ae_2=\frac{\lambda_{12}}{\mu_1}e_2$ being a valid automorphism. A similar computation yields $(\Delta_{1,2},0)+(0,\omega_2)\cong ({\Delta}_{1,2}',0)$.

For the general case, we show that $(\Delta_{1,2},0)+(0,\omega_1)+(0,\omega_2)\cong (\Delta_{1,2}',0)$. However, for the latter cocycle, we explicitly set $\Delta_{1,2}'(e_1,e_2)=1$ for our convenience. We may do so as $(\Delta_{1,2},0)$ (where $\Delta_{1,2}=\lambda_{12}\in \K\setminus\{0\}$) belongs to one single equivalence class. Indeed, we find that equating the conditions in \eqref{conditions} reduces to
\begin{equation*}
b_2^2\mu_1+a_2^2\mu_2=a_2b_2\lambda_{12}.
\end{equation*}
Setting $d=1$, we obtain a solution for $b_2$ by invoking \textit{algebraic closure}. We verify that the values obtained do not yield det($A$)$=0$, and hence we have a valid automorphism.

\vspace{0.5 em}

The second case is trivial, as $(0,0)$ is the only cocycle. For the the third case $L=\mathbf{L^3_{1|1}}$, the general form of the automorphism is
\begin{equation*}
Ae_1=a_1e_1, ~Ae_2=b_2e_2.
\end{equation*}

We have dim($H^2(L;\K)$)$=2$, with the basis being given by $\{(0,\omega_2) ~(even),(\Delta_{1,2},0) ~(odd)\}$. We verify that $(0,\omega_2) \cong (0,\omega_2')$ with the automorphism $Ae_1=e_1, ~Ae_2=\sqrt{\frac{\mu_2}{\mu_2'}}e_2$. Similarly, we check $(\Delta_{1,2},0)\cong (\Delta_{1,2}',0)$, where $Ae_1=e_1, ~Ae_2=\frac{\lambda_{12}}{\lambda_{12}'}e_2$ is a valid automorphism. 

\vspace{0.5 em}

For the fourth and final case $L=\mathbf{L^2_{2|0}}$, the general form of the automorphism is given by
\begin{equation*}
Ae_1=a_1e_1+a_2e_2, ~Ae_2=b_1e_1+b_2e_2.
\end{equation*}

We find that we only have $(\Delta_{1,2},0)$ as a basis element. It is even. We check that $(\Delta_{1,2},0) \cong (\Delta_{1,2}',0)$, using the automorphism $Ae_1=e_1, ~Ae_2=\frac{\lambda_{12}}{\lambda_{12}'}e_2$.
\end{proof}

\subsection{Building the Extensions.}\label{4.3} As elucidated in Theorem \ref{indthm}, we use the central extensions method with non-equivalent scalar 2-cocycles to produce a complete list of nilpotent Lie superalgebras of dimension 3, though the method yields isomorphic cases which shall be dealt with in the succeeding subsection. For each extension, we tabulate the scalar 2-cocycle used, the parity of the added element $X$, as well as the new Lie bracket and squaring. In alignment with Step 2 of Theorem \ref{thm1}, and Theorem \ref{indthm}, the relations on the new superalgebras are derived as follows:
\begin{align}
  [x,y]_{\mathrm{New}}&=[x,y]_{\mathrm{Old}}+
  \Delta(x,y) X, ~\forall x,y\in L\\
  s_{\mathrm{New}}(x)&=s_{\mathrm{Old}}(x)+\omega(x)X, ~\forall x\in L_{\od} .
\end{align}

We also recall from Theorem \ref{indthm}, that the extension of a Lie superalgebra $L$ by the ideal generated by the added element $X$ may be viewed as an extension of $L$ by the underlying field using a cocycle of the same parity as $X$. Hence, we always have $|X|=|(\Delta,\omega)|=|\Delta|$ (see \eqref{cograd}).

\vspace{0.5 em}

{\renewcommand{\arraystretch}{1.5} 
\small
\begin{center}
    \makebox[\textwidth][c]{
    \begin{tabular}{|c|c|c|c|c|c|}
    \hline
          \textbf{Name} &\textbf{sdim} & \textbf{Cocycle}& \textbf{Added element} &\textbf{Lie Bracket} &\textbf{Squaring} \\\hline
          ${\bf L^a_{0|3}}$ & $(0|3)$ & $0$ & $X$ odd & $[\cdot,\cdot]=0$ & $s(\cdot)=0$ \\\hline
         ${\bf L^a_{1|2}}$ &$(1|2)$ & $0$ & $X$ even & $[\cdot,\cdot]=0$ & $s(\cdot)=0$  \\\hline
         ${\bf L^b_{1|2}}$& $(1|2)$  & $(0,\omega_1)$   & $X$ even & $[\cdot,\cdot]=0$ & $s(e_1)=X$ \\\hline
         ${\bf L^c_{1|2}}$& $(1|2)$  & $(\Delta_{1,2},0)$   & $X$ even & $[e_1,e_2]=X$ & $s(\cdot)=0$\\\hline
    \end{tabular}
    }
   \\~\\ Lie superalgebras obtained by extensions of ${\bf L^1_{0|2}}$.
\end{center}

\vspace{0.5 em}

\begin{center}
    \makebox[\textwidth][c]{
    \begin{tabular}{|c|c|c|c|c|c|}
    \hline
          \textbf{Name} &\textbf{sdim} & \textbf{Cocycle}& \textbf{Added element} &\textbf{Lie Bracket} &\textbf{Squaring} \\\hline
          ${\bf L^a_{2|1}}$ & $(2|1)$ & $0$ & $X$ even & $[\cdot,\cdot]=0$ & $s(e_2)=e_1$ \\\hline
         ${\bf L^d_{1|2}}$ &$(1|2)$ & $0$ & $X$ odd & $[\cdot,\cdot]=0$ & $s(e_2)=e_1$  \\\hline
    \end{tabular}
    }
   \\~\\ Lie superalgebras obtained by extensions of ${\bf L^2_{1|1}}$.
\end{center}

\vspace{0.5 em}

\begin{center}
    \makebox[\textwidth][c]{
    \begin{tabular}{|c|c|c|c|c|c|}
    \hline
          \textbf{Name} &\textbf{sdim} & \textbf{Cocycle}& \textbf{Added element} &\textbf{Lie Bracket} &\textbf{Squaring} \\\hline
          ${\bf L^b_{2|1}}$ & $(2|1)$ & $0$ & $X$ even & $[\cdot,\cdot]=0$ & $s(\cdot)=0$ \\\hline
         ${\bf L^e_{1|2}}$ &$(1|2)$ & $0$ & $X$ odd & $[\cdot,\cdot]=0$ & $s(\cdot)=0$  \\\hline
         ${\bf L^c_{2|1}}$& $(2|1)$  & $(0,\omega_2)$   & $X$ even & $[\cdot,\cdot]=0$ & $s(e_2)=X$ \\\hline
         ${\bf L^f_{1|2}}$& $(1|2)$  & $(\Delta_{1,2},0)$   & $X$ odd & $[e_1,e_2]=X$ & $s(\cdot)=0$\\\hline
    \end{tabular}
    }
   \\~\\ Lie superalgebras obtained by extensions of ${\bf L^3_{1|1}}$.
\end{center}

\vspace{0.5 em}

\begin{center}
    \makebox[\textwidth][c]{
    \begin{tabular}{|c|c|c|c|c|c|}
    \hline
          \textbf{Name} &\textbf{sdim} & \textbf{Cocycle}& \textbf{Added element} &\textbf{Lie Bracket} &\textbf{Squaring} \\\hline
          ${\bf L^a_{3|0}}$ & $(3|0)$ & $0$ & $X$ even & $[\cdot,\cdot]=0$ & $s(\cdot)=0$ \\\hline
         ${\bf L^d_{2|1}}$ &$(2|1)$ & $0$ & $X$ odd & $[\cdot,\cdot]=0$ & $s(\cdot)=0$  \\\hline
         ${\bf L^b_{3|0}}$& $(3|0)$  & $(\Delta_{1,2},0)$   & $X$ even & $[e_1,e_2]=X$ & $s(\cdot)=0$ \\\hline
    \end{tabular}
    }
   \\~\\ Lie superalgebras obtained by extensions of ${\bf L^2_{2|0}}$.
\end{center}
\normalsize
\vspace{0.5 em}

\subsubsection{Removing Isomorphic Cases}\label{iso3} We find that a number of Lie superalgebras derived in Section \ref{4.3} are isomorphic. For a complete, \emph{non-redundant} list, duplicates must be removed. We find the following isomorphic cases up to a simple relabeling of the basis elements.

\vspace{0.5 em}

\noindent$\underline{\sdim(L)=(1|2)}$: $\mathbf{L^{a}_{1|2}}   \cong$ $\mathbf{L^{e}_{1|2}}$, $\mathbf{L^{b}_{1|2}} \cong$ $\mathbf{L^{d}_{1|2}}$

\vspace{0.5 em}

\noindent$\underline{\sdim(L)=(2|1)}$: $\mathbf{L^{a}_{2|1}}  \cong$ $\mathbf{L^{c}_{2|1}} $, $\mathbf{L^{b}_{2|1}} \cong$ $\mathbf{L^{d}_{2|1}}$

\vspace{0.5 em}

As for the rest, we determine that they are non-isomorphic by explicitly yielding contradictory isomorphism conditions. We omit these computations, though they are easily verifible. Moreover, most of these can be deduced from their second cohomolgies as listed in \ref{5.0.1} which is an invariant. Hence, we list the \emph{full, non-redundant} list in dimension 3 as follows.

\sssbegin{Theorem}\label{dim3}
The classification of 3-dimensional nilpotent Lie superalgebras over an algebraically closed field of characteristic 2 is given by:

\vspace{0.5 em}

\noindent$\underline{\sdim(L)=(0|3)}$:

    $\mathbf{L^{1}_{0|3}} = \langle 0|e_{1},e_{2},e_3\rangle \ (abelian)$.\\
    
\noindent$\underline{\sdim(L)=(1|2)}$:

    $\mathbf{L^{1}_{1|2}} = \langle e_{1}|e_{2},e_{3}\rangle \ (abelian)$;
    
    $\mathbf{L^{2}_{1|2}} = \langle e_1|e_{2},e_{3}; s(e_2)=e_1 \rangle$;
    
    $\mathbf{L^{3}_{1|2}} = \langle e_1|e_{2},e_{3}; [e_2,e_3]=e_1 \rangle$;
    
    $\mathbf{L^{4}_{1|2}} = \langle e_{1}|e_{2},e_3; [e_1,e_2]=e_3\rangle$.\\

\noindent$\underline{\sdim(L)=(2|1)}$:

    $\mathbf{L^{1}_{2|1}} = \langle e_{1},e_{2}|e_3\rangle \ (abelian)$.
    
    $\mathbf{L^{2}_{2|1}} = \langle e_{1},e_{2}|e_{3}; s(e_{3}) = e_{1}\rangle$.\\

\noindent$\underline{\sdim(L)=(3|0)}$:

    $\mathbf{L^{1}_{3|0}} = \langle e_{1},e_{2},e_3|0\rangle \ (abelian)$;
    
    $\mathbf{L^{2}_{3|0}} = \langle e_{1},e_{2},e_3|0; [e_1,e_2]=e_3\rangle$.
\end{Theorem}

\section{Classification in Dimension 4} \label{class4}

We now adapt the method applied in Section 4, to the the list of  superalgebras obtained in Theorem \ref{dim3}, to classify the nilpotent Lie superalgebras in dimension 4. Again, utilizing \texttt{SuperLie} (see \cite{Kr}), we compute the following 2-cocyle list. 

\sssbegin{Lemma}\label{5.0.1} Suppose that L is a nilpotent Lie superalgebra with $\dim(L)=3$ over an algebraically closed field of characteristic 2. Then, the scalar 2-cocyle basis of L is given by
\end{Lemma}

\vspace{0.5 em}
\small
\begin{center}
    \makebox[\textwidth][c]{
    \begin{tabular}{|c|c|c|c|}
    \hline
          \textbf{Name} &\textbf{sdim} & \textbf{Even Cocycles}& \textbf{Odd Cocycles} \\\hline
          $\mathbf{L^{1}_{0|3}}$  & $(0|3)$ & $(\Delta_{1,2},0),(\Delta_{1,3},0),(\Delta_{2,3},0),(0,\omega_1),(0,\omega_2),(0,\omega_3)$ & --- \\\hline
          $\mathbf{L^{1}_{1|2}}$  & $(1|2)$ & $(\Delta_{2,3},0)  ,(0,\omega_2)  ,(0,\omega_3)$ & $(\Delta_{1,2},0),(\Delta_{1,3},0)$ \\\hline
          $\mathbf{L^{2}_{1|2}}$  & $(1|2)$ & $(\Delta_{2,3},0), (0,\omega_3)$ & --- \\\hline
          $\mathbf{L^{3}_{1|2}}$  & $(1|2)$ & $ (0,\omega_2), (0,\omega_3)$ & --- \\\hline
          $\mathbf{L^{4}_{1|2}}$  & $(1|2)$ & $ (0,\omega_2), (0,\omega_3)$ & $(\Delta_{1,3},0)$ \\\hline
          $\mathbf{L^{1}_{2|1}}$  & $(2|1)$ & $ (\Delta_{1,2},0), (0,\omega_3)$ & $(\Delta_{1,3},0), (\Delta_{2,3},0)$ \\\hline
          $\mathbf{L^{2}_{2|1}}$  & $(2|1)$ & --- & $(\Delta_{2,3},0)$ \\\hline
          $\mathbf{L^{1}_{3|0}}$  & $(3|0)$ & $(\Delta_{1,2},0)$,$(\Delta_{1,3},0)$,$(\Delta_{2,3},0)$ & --- \\\hline
          $\mathbf{L^{2}_{3|0}}$  & $(3|0)$ & $(\Delta_{1,3},0)$,$(\Delta_{2,3},0)$ & --- \\\hline
          
    \end{tabular}
    }
   \\~\\ Scalar $2$-cocyle basis of nilpotent Lie superalgebras of dimension $3$.
\end{center}
\normalsize
\vspace{0.5 em}

\subsection{Automorphism Computation for \texorpdfstring{$\mathrm{dim}(L)=3$}{dim(L)=3}.}

An analogue to \ref{4.2.1}, we compute the equivalence classes of the scalar coycles for nilpotent Lie superalgebras of dimension 3.

\sssbegin{Lemma}\label{5.1.1}
Suppose that $L$ is a nilpotent Lie superalgebra with $\sdim(L)=(0|3)$ over an algebraically closed field of characteristic $p=2$. The equivalence classes of scalar $2$-cocycles on $L$ are given by: $(0,0),~(\Delta_{1,2},0), ~(\Delta_{1,2},\omega_3), ~(0,\omega_1)$.
\end{Lemma}

\begin{proof}
The general form of the automorphism $A \in \Aut(L)$ is given by
\begin{equation}\label{autform}Ae_1=a_1e_1+a_2e_2+a_3e_3,~Ae_2=b_1e_1+b_2 e_2+b_3 e_3,~Ae_3=c_1e_1+c_2e_2+c_3e_3,\end{equation}
where $a_i,b_i,c_i \in \K$ for $i=1,2,3$. There are no additional conditions on $A$ except $\det(A)\not =0$. Moreover, we find that $\dim\left(H^2(L;\K)\right)=6$ with the basis being given by
$$\left\lbrace(\Delta_{1,2},0),~(\Delta_{1,3},0),~(\Delta_{2,3},0),~(0,\omega_1),~(0,\omega_2),~(0,\omega_3)\right\rbrace.$$
Using the same method as detailed in Section \ref{method}, we very quickly find that $3$ distinct equivalence classes arise: $(\Delta_{1,2},0) \cong (\Delta_{1,3},0) \cong (\Delta_{2,3},0)$, $~(0,\omega_1) \cong (0,\omega_2) \cong (0,\omega_3)$, and $(\Delta_{1,2},\omega_3) \cong (\Delta_{1,3},\omega_2) \cong (\Delta_{2,3},\omega_1)$. Then, we systematically prove that every combination of the basis elements must belong to one of these three classes. We find that the computations are remarkably lengthy and contain many technicalities. Such was the case in \cite[Lemma 5.3.1]{BE}, and one imagines that it is so more generally for adaptations of the Skjelbred-Sund method for certain abelian cases. As is the practice, we omit the computations entirely. However, we still sketch a rough outline for the most general case, which contains a number of techniques that are pertinent to the other cases as well as some novel ones which are unique to itself.

We find that setting up the conditions in \eqref{conditions} gives way to a set of parameters on \eqref{autform} such that $(\varphi, \omega) = (\Delta_{1,2},0)+(\Delta_{1,3},0)+(\Delta_{2,3},0)+(0,\omega_1)+(0,\omega_2)+(0,\omega_3) \cong (\Delta_{1,2}', \omega_3')$ if and only if 
\begin{equation}\label{relation}
\lambda_{12}^2\mu_{3}+\lambda_{13}^2\mu_{2}+\lambda_{23}^2\mu_{1} = \lambda_{12}\lambda_{13}\lambda_{23}
\end{equation}
does not hold. In other words, all cocycles of the general form that do not satisfy the relation in \eqref{relation} are automorphic to $(\Delta_{1,2}, \omega_3)$. 

Now, for the ones that do satisfy \eqref{relation}, we find that it is difficult to find a direct automorphism to $(\Delta_{1,2},0)$. Hence, we first show that that all such cocycles are automorphic to the general case where each of the basis elements is set to 1, i.e. $\lambda_{ij}=\mu_k=1$ where $i \neq j$ and $i,j,k \in \{1,2,3\}$. We clearly see that the latter satisfies \eqref{relation}. The conditions in \eqref{conditions} yield a daunting system of equations. By carefully setting certain parameters, and using \eqref{relation}, we relax precise conditions and find a valid automorphism.

Finally, we attempt to find an automorphism between this latter cocyle and $(\Delta_{1,2},0)$. We again set up the conditions, and find that they reduce to 
$$
\frac{a_1}{b_1}+\frac{b_1}{a_1}=1,\quad 
\frac{a_1}{c_1}+\frac{c_1}{a_1}=1,\quad 
\frac{b_1}{c_1}+\frac{c_1}{b_1}=1.
$$

Setting $\frac{a_1}{b_1}=x$, $\frac{a_1}{c_1}=y$, and $\frac{b_1}{c_1}=\frac{y}{x}$, we yield the following set of equations in $x$ and $y$:
$$x^2+x+1=0, ~y^2+y+1=0, ~\frac{y}{x}+\frac{x}{y}+1=0.$$
Using algebraic closure on the equation in $x$, we yield a solution $\alpha \neq 0,1$. We see that $1+\alpha$ also has to satisfy it, and hence, we set $x=\alpha$ and $y=1+\alpha$ which is a valid solution for the set of equations. From there it is easy to set the values of the parameters and determine $\det(A)\neq0$.\end{proof}

\sssbegin{Lemma}
Suppose that $L$ is a nilpotent Lie superalgebra with $\sdim(L)=(1|2)$ over an algebraically closed field of characteristic $p=2$. The equivalence classes of scalar $2$-cocycles on $L$ are given by

\begin{itemize}
    \item[] $L=\mathbf{L^{1}_{1|2}}$: $(0,0), ~(0,\omega_2), ~(\Delta_{2,3},0), ~(\Delta_{1,2},0)$;
    
    \item[] $L=\mathbf{L^{2}_{1|2}}$: $(0,0), ~(0,\omega_3), ~(\Delta_{2,3},0), ~(\Delta_{2,3},\omega_3)$;
    
    \item[] $L=\mathbf{L^{3}_{1|2}}$: $(0,0), ~(0,\omega_2), ~(0,\omega_2)+(0,\omega_3)$; 

    \item[] $L=\mathbf{L^{4}_{1|2}}$: $(0,0), ~(0,\omega_2), ~(0,\omega_3), ~(\Delta_{1,3},0)$.
\end{itemize}
\end{Lemma}
\begin{proof}
The general form of the automorphism $A \in \Aut(L)$ is given by,
$$Ae_1=a_1e_1,~Ae_2=b_2 e_2+b_3 e_3,~Ae_3=c_2e_2+c_3e_3,\text{ where }a_1,b_2,b_3,c_2,c_3\in\K.$$
The condition on $\det(A)$ is given by $a_1(b_2c_3+c_2b_3)\neq0$.
\begin{itemize}
\item $L=\mathbf{L^{1}_{1|2}}$. No further conditions of $A$. We have $\dim\left(H^2_{}(L;\K)\right)=5$, with two odd cocycles, $(\Delta_{1,2},0),(\Delta_{1,3},0)$, and three even cocycles, $(\Delta_{2,3},0),(0,\omega_2),(0,\omega_3)$. One verifies $(\Delta_{1,2},0) \cong (\Delta_{1,3},0)$, and $(\Delta_{1,2},0)+(\Delta_{1,3},0) \cong (\Delta_{1,2}',0)$, the latter given by the automorphism,
$$Ae_1=\frac{1}{\lambda_{12}'}e_1, ~Ae_2=\lambda_{12}e_2, ~Ae_3=\lambda_{13}e_2+e_3.$$
Hence, there is only one non-trivial equivalence class of odd cocycles: $(\Delta_{1,2},0)$. For the even cocycles, one verifies the following: $(0,\omega_2)\cong (0,\omega_3')$, $(\Delta_{2,3}0) \ncong (0,\omega_2)$,  $(\Delta_{2,3},0)+(0,\omega_2)\cong (\Delta_{2,3}',0)$, $(\Delta_{2,3},0)+(0,\omega_3)\cong (\Delta_{2,3}',0)$. Moreover, we find $(0,\omega_2)+(0,\omega_3) \cong (0,\omega_2')$, with the automorphism
$$Ae_1=e_1, ~Ae_2= \sqrt{\frac{\mu_2}{\mu_2'}}e_2, ~Ae_3=\sqrt{\frac{\mu_3}{\mu_2'}}e_2+e_3.$$ Again, we set $\Delta_{2,3}'(e_2,e_3)=1$ as all $(\Delta_{2,3},0)$ belong to the same equivalence class, and verify $(\Delta_{2,3},0)+(0,\omega_2)+(0,\omega_3)\cong (\Delta_{2,3}',0)$, given by the automorphism
$$Ae_1=e_1, ~Ae_2=e_2+\mu_2e_3, ~Ae_3=c_2e_2+\frac{\mu_3}{c_2}e_3,$$
where $c_2$ is a solution of the equation $\mu_2c_2^2+\lambda_{23}c_2+\mu_3=0$. Hence, there are two non-trivial equivalence classes of even cocycles: $(\Delta_{1,2},0)$ and $(0,\omega_2)$.

\vspace{0.5 em}

\item $L=\mathbf{L^{2}_{1|2}}$. The additional conditions on $A$ are $a_1=b_2^2$ and $c_2=0$. We have $\dim\left(H^2_{}(L;\K)\right)=2$, with $(\Delta_{2,3},0)$ and $(0, \omega_3)$ both even cocycles. It is easily checked that $(\Delta_{2,3},0) \ncong (0, \omega_3)$ and $(\Delta_{2,3},0) \ncong (\Delta_{2,3},0) +(0, \omega_3)$. One may also verify that $(\Delta_{2,3},0)\cong(\Delta_{2,3}',0)$ and $(0,\omega_3) \cong (0,\omega_3')$. We finally check that $(\Delta_{2,3},\omega_3)\cong(\Delta_{2,3}',\omega_3')$, which is given by the automorphism
$$Ae_1=\frac{\lambda_{23}^2\mu_3'}{(\lambda_{23}')^2\mu_3}e_1, ~Ae_2=\sqrt{\frac{\lambda_{23}^2\mu_3'}{(\lambda_{23}')^2\mu_3}}e_2, ~Ae_3=\sqrt{\frac{\mu_3}{\mu_3'}}e_3.$$ 
Hence, there are three non-trivial equivalence classes of even cocycles: $(0,\omega_3), (\Delta_{2,3},0)$ and $(\Delta_{2,3},\omega_3)$.

\vspace{0.5 em}

\item $L=\mathbf{L^{3}_{1|2}}$. The additional conditions on $A$ are $a_1=b_2c_3+b_3c_2$, $~b_2b_3=0$, and $c_2c_3=0$. We have $\dim\left(H^2_{}(L;\K)\right)=2$ with $(0,\omega_2)$ and $(0, \omega_3)$ being the only two cocyles, both of which are even. It can be checked that $(0,\omega_2) \cong (0,\omega_3')$, but $(0,\omega_2) \ncong (0,\omega_2)+(0,\omega_3)$. However, we find that $(0,\omega_2)+(0,\omega_3) \cong (0,\omega_2')+(0,\omega_3')$, with the automorphism
$$Ae_1= \sqrt{\frac{\mu_2\mu_3}{\mu_2' \mu_3'}}e_1, ~Ae_2=\sqrt{\frac{\mu_2}{\mu_2'}}e_2, ~Ae_3=\sqrt{\frac{\mu_3}{\mu_3'}}e_3.$$ 
Hence, there are two non-trivial equivalence classes of even cocycles: $(0,\omega_2)$ and $(0,\omega_2)+(0,\omega_3)$.

\vspace{0.5 em}

\item $L=\mathbf{L^{4}_{1|2}}$. The additional conditions on $A$ are $c_2=0$ and $a_1b_2=c_3$. We have $\dim\left(H^2_{}(L;\K)\right)=3$, with two even cocycles $(0,\omega_2)$ and $(0,\omega_3)$, and one odd cocycle, $(\Delta_{1,3},0)$. For the odd case, it is easily checked that $(\Delta_{1,3},0)\cong(\Delta_{1,3}',0)$, and that there is only one non-trivial equivalence class for odd cocycles, that being $(\Delta_{1,3},0)$. For the even case, due to the additional conditions imposed on $A$, we find that $(0,\omega_2)\ncong(0,\omega_3)$. However, we verify $(0,\omega_2) \cong (0,\omega_2')$ and $(0,\omega_3) \cong (0,\omega_3')$. Lastly, we find that $  (0,\omega_2)+(0,\omega_3) \cong (0,\omega_3')$, under the automorphism
$$Ae_1=e_1, ~Ae_2=\sqrt{\frac{\mu_3}{\mu_3'}}e_2+\sqrt{\frac{\mu_2}{\mu_3'}}e_3, ~Ae_3=\sqrt{\frac{\mu_3}{\mu_3'}}e_3.$$
Hence, there are two non-trivial equivalence classes for even cocyles: $(0,\omega_2)$ and $(0,\omega_3)$. \qed
\end{itemize}\noqed \end{proof}

\sssbegin{Lemma}
Suppose that $L$ is a nilpotent Lie superalgebra with $\sdim(L)=(2|1)$ over an algebraically closed field of characteristic $p=2$. The equivalence classes of scalar $2$-cocycles on $L$ are given by

\begin{itemize}
    \item[]  $\mathbf{L^{1}_{2|1}}$: $(0,0), ~(\Delta_{1,2},0), ~(0,\omega_3), ~(\Delta_{1,2},\omega_3), ~(\Delta_{1,3},0)$; 
    
    \item[] $\mathbf{L^{2}_{2|1}}$: $(0,0), ~(\Delta_{2,3},0).$
\end{itemize}
\end{Lemma}

\begin{proof}
The general form of the automorphism $A \in \Aut(L)$ is given by,
$$Ae_1=a_1e_1+a_2e_2,~Ae_2=b_1 e_1+b_2 e_2,~Ae_3=c_3e_3,\text{ where }a_1,a_2,b_1,b_2,c_3\in\K.$$
Clearly, the condition on $\det(A)$ is given by $c_3(a_1b_2+b_1a_2)\neq0$.
\begin{itemize}
\item $L=\mathbf{L^1_{2|1}}$. No further conditions on $A$. We have $\dim\left(H^2_{}(L;\K)\right)=4$, with two even cocycles,  $(\Delta_{1,2},0), (0,\omega_3)$, and two odd cocycles, $(\Delta_{1,3},0)$, $(\Delta_{2,3},0)$. One easily verifies that $(\Delta_{1,2},0) \ncong (0, \omega_3)$, $(0,\omega_3) \cong (0,\omega_3')$, and that $(\Delta_{1,2},\omega_3) \ncong (\Delta_{1,2}',0)$. We find that $(\Delta_{1,2},\omega_3) \cong (\Delta_{1,2}',\omega_3')$, given by the automorphism
$$Ae_1=\lambda_{12}e_1, ~Ae_2=\frac{1}{\lambda_{12}'}e_2, ~Ae_3=\sqrt{\frac{\mu_3}{\mu_3'}}e_3.$$
Therefore, we have three non-trivial equivalence classes of even cocycles: $(0,\omega_3)$, $(\Delta_{1,2},0)$ and $(\Delta_{1,2}, \omega_3)$. For the odd cocycles, we have $(\Delta_{2,3},0) \cong (\Delta_{1,3}',0)  $, given by 
$$Ae_1=e_2, ~Ae_2=\lambda_{23}e_1, ~Ae_3=\frac{1}{\lambda_{13}'}e_3.$$
Moreover, we find check that $(\Delta_{1,3},0)+(\Delta_{2,3},0) \cong (\Delta_{1,3}',0)$ by the automorphism,
$$Ae_1=\lambda_{13}e_1, ~Ae_2=\lambda_{23}e_1+e_2, ~Ae_3=\frac{1}{\lambda_{13}'}e_3.$$
Hence, we have one non-trivial equivalence class of odd cocycles: $(\Delta_{1,3},0)$. 

\item $L=\mathbf{L^{2}_{2|1}}$. The additional conditions on $A$ are $a_2=0$ and $c_3^{2}=a_1$. We only have one non-trivial cocycle: $(\Delta_{2,3},0)$ which is odd. We check that $(\Delta_{2,3},0) \cong (\Delta_{2,3}',0)$, given by the automorphism
$$Ae_1=\frac{1}{({\lambda_{23}'})^2}e_1, ~Ae_2=\lambda_{23}e_2, ~Ae_3=\frac{1}{\lambda_{23}'}e_3.$$
Hence, we have one non-trivial equivalence class of odd cocycles: $(\Delta_{2,3},0)$. \qed
\end{itemize}\noqed \end{proof}

\sssbegin{Lemma}
Suppose that $L$ is a nilpotent Lie superalgebra with $\sdim(L)=(3|0)$ over an algebraically closed field of characteristic $p=2$. The equivalence classes of scalar $2$-cocycles on $L$ are given by

\begin{itemize}
    \item[] $\mathbf{L^{1}_{3|0}}$: $(0,0), ~(\Delta_{1,2},0)$;
    
    \item[] $\mathbf{L^{2}_{3|0}}$: $(0,0), ~(\Delta_{1,3},0)$. 
\end{itemize}
\end{Lemma}

\begin{proof}
The general form of an automorphism $A\in \Aut(L)$ is given by
$$Ae_1=a_1e_1+a_2e_2+a_3e_3,~Ae_2=b_1e_1+b_2 e_2+b_3e_3,~Ae_3=c_1e_1+c_2e_2+c_3e_3,$$
$\text{ where }a_i,b_i,c_i\in\K.$ for $i=1,2,3$. 

\begin{itemize}

\item $L=\mathbf{L^{1}_{3|0}}$. No conditions on $A$ except $\det (A)\neq 0$. We have $\dim\left(H^2_{}(L;\K)\right)=3$ with the basis elements $(\Delta_{1,2},0)$, $(\Delta_{1,3},0)$, $(\Delta_{2,3},0)$ which are all even. We get that all linear combinations of these basis elements are equivalent to $(\Delta_{1,2},0)$, except for $(\Delta_{1,2},0)+(\Delta_{1,3},0)+(\Delta_{2,3},0)$, as a consequence of Lemma \ref{5.1.1}.  Indeed, we find the same conditions arise minus those imposed by the squaring. We set $\Delta_{1,2}'=1$. We find $(\Delta_{1,2},0)+(\Delta_{1,3},0)+(\Delta_{2,3},0) \cong (\Delta_{1,2}',0)  $ independently using the automorphism
$$Ae_1=e_2, ~Ae_2= \lambda_{12}e_1+ \left(\frac{\lambda_{23}+1}{\lambda_{13}}\right)e_2, ~Ae_3=\frac{1}{\lambda_{12}}e_1+\lambda_{13}e_2+e_3.$$
Hence, $(\Delta_{1,2},0)$ is the only non-trivial equivalence class of scalar cocycles of $L$.

\vspace{0.5 em}

\item $L=\mathbf{L^{2}_{3|0}}$. The conditions on $A$ are $c_1,c_2=0$ and $c_3=a_1b_2+b_1a_2$. We have $\dim\left(H^2_{}(L;\K)\right)=2$ with the basis elements $(\Delta_{1,3},0)$, $(\Delta_{2,3},0)$ which are both even. We have $(\Delta_{1,3},0) \cong (\Delta_{2,3}',0)$, followed by $(\Delta_{1,3},0)+(\Delta_{2,3},0) \cong (\Delta_{1,3}',0)$, the latter being given by the automorphism:

$$Ae_1=\lambda_{13}e_1, ~Ae_2=\lambda_{23}e_1+\frac{1}{\lambda_{13}\lambda_{13}'}e_2, ~Ae_3=\frac{1}{\lambda_{13}'}e_3.$$ 

Hence, $(\Delta_{1,3},0)$ is the only non-trivial equivalence class of scalar cocycles of $L$. \qed  
\end{itemize} \noqed
\end{proof}

\subsection{Buidling the Extensions for \texorpdfstring{$\mathrm{dim}(L)=4$.}{dim(L)=4}}\label{5.2} As in the presentation of Section \ref{4.3}, we present tables of the extensions of $nilpotent$ Lie superalgebras in Theorem \ref{dim3} to produce a complete list of nilpotent Lie superalgebras in dimension 4, though isomorphic cases may arise, which are dealt with in the next subsection.\\
\small
{\renewcommand{\arraystretch}{1.5}

\begin{center}
    \makebox[\textwidth][c]{
    \begin{tabular}{|c|c|c|c|c|c|}
    \hline
          \textbf{Name} &\textbf{sdim} & \textbf{Cocycle}& \textbf{Added element} &\textbf{Lie Bracket} &\textbf{Squaring} \\\hline
          ${\bf L^a_{0|4}}$ & $(0|4)$ & $0$ & $X$ odd & $[\cdot,\cdot]=0$ & $s(\cdot)=0$ \\\hline
         ${\bf L^a_{1|3}}$ &$(1|3)$ & $0$ & $X$ even & $[\cdot,\cdot]=0$ & $s(\cdot)=0$  \\\hline
         ${\bf L^b_{1|3}}$& $(1|3)$  & $(0,\omega_1)$   & $X$ even & $[\cdot,\cdot]=0$ & $s(e_1)=X$ \\\hline
         ${\bf L^c_{1|3}}$& $(1|3)$  & $(\Delta_{1,2},0)$   & $X$ even & $[e_1,e_2]=X$ & $s(\cdot)=0$\\\hline
         ${\bf L^d_{1|3}}$& $(1|3)$  & $(\Delta_{1,2},\omega_3)$   & $X$ even & $[e_1,e_2]=X$ & $s(e_3)=X$\\\hline
    \end{tabular}
    }
   \\~\\ Lie superalgebras obtained by extensions of ${\bf L^1_{0|3}}$.
\end{center}
}

\vspace{0.5 em}

{\renewcommand{\arraystretch}{1.5}
\begin{center}
    \makebox[\textwidth][c]{
    \begin{tabular}{|c|c|c|c|c|c|}
    \hline
          \textbf{Name} &\textbf{sdim} & \textbf{Cocycle}& \textbf{Added element} &\textbf{Lie Bracket} &\textbf{Squaring} \\\hline
          ${\bf L^a_{4|0}}$ & $(4|0)$ & $0$ & $X$ even & $[\cdot,\cdot]=0$ & $s(\cdot)=0$ \\\hline
         ${\bf L^a_{3|1}}$ &$(3|1)$ & $0$ & $X$ odd & $[\cdot,\cdot]=0$ & $s(\cdot)=0$  \\\hline
         ${\bf L^b_{4|0}}$& $(4|0)$  & $(\Delta_{1,2},0)$   & $X$ even & $[e_1,e_2]=X$ & $s(\cdot)=0$ \\\hline
    \end{tabular}
    }
   \\~\\ Lie superalgebras obtained by extensions of ${\bf L^1_{3|0}}$.
\end{center}
}

\vspace{0.5 em}

{\renewcommand{\arraystretch}{1.5}
\begin{center}
    \makebox[\textwidth][c]{
    \begin{tabular}{|c|c|c|c|c|c|}
    \hline
          \textbf{Name} &\textbf{sdim} & \textbf{Cocycle}& \textbf{Added element} &\textbf{Lie Bracket} &\textbf{Squaring} \\\hline
          ${\bf L^c_{4|0}}$ & $(4|0)$ & $0$ & $X$ even & $[e_1,e_2]=e_3$ & $s(\cdot)=0$ \\\hline
         ${\bf L^b_{3|1}}$ & $(3|1)$ & $0$ & $X$ odd & $[e_1,e_2]=e_3$ & $s(\cdot)=0$ \\\hline
         ${\bf L^d_{4|0}}$& $(4|0)$  & $(\Delta_{1,3},0)$   & $X$ even & $[e_1,e_2]=e_3$, $~[e_1,e_3]=X$ & $s(\cdot)=0$ \\\hline
    \end{tabular}
    }
   \\~\\ Lie superalgebras obtained by extensions of ${\bf L^2_{3|0}}$.
\end{center}
}
\vspace{0.5 em}

{\renewcommand{\arraystretch}{1.5}
\begin{center}
    \makebox[\textwidth][c]{
    \begin{tabular}{|c|c|c|c|c|c|}
    \hline
          \textbf{Name} &\textbf{sdim} & \textbf{Cocycle}& \textbf{Added element} &\textbf{Lie Bracket} &\textbf{Squaring} \\\hline
          ${\bf L^c_{3|1}}$ & $(3|1)$ & $0$ & $X$ even & $[\cdot,\cdot]=0$ & $s(\cdot)=0$ \\\hline
          ${\bf L^a_{2|2}}$ & $(2|2)$ & $0$ & $X$ odd & $[\cdot,\cdot]=0$ & $s(\cdot)=0$ \\\hline
          ${\bf L^d_{3|1}}$ & $(3|1)$ & $(\Delta_{1,2},0)$ & $X$ even & $[e_1,e_2]=X$ & $s(\cdot)=0$ \\\hline
          ${\bf L^e_{3|1}}$ & $(3|1)$ & $(0,\omega_3)$ & $X$ even & $[\cdot,\cdot]=0$ & $s(e_3)=X$ \\\hline
          ${\bf L^b_{2|2}}$ & $(2|2)$ & $(\Delta_{1,3},0)$ & $X$ odd & $[e_1,e_3]=X$ & $s(\cdot)=0$ \\\hline
          ${\bf L^f_{3|1}}$ & $(3|1)$ & $(\Delta_{1,2},\omega_3)$ & $X$ even & $[e_1,e_2]=X$ & $s(e_3)=X$ \\\hline
    \end{tabular}
    }
   \\~\\ Lie superalgebras obtained by extensions of ${\bf L^1_{2|1}}$.
\end{center}
}
\vspace{0.5 em}

{\renewcommand{\arraystretch}{1.5}
\begin{center}
    \makebox[\textwidth][c]{
    \begin{tabular}{|c|c|c|c|c|c|}
    \hline
          \textbf{Name} &\textbf{sdim} & \textbf{Cocycle}& \textbf{Added element} &\textbf{Lie Bracket} &\textbf{Squaring} \\\hline
          ${\bf L^g_{3|1}}$ & $(3|1)$ & $0$ & $X$ even & $[\cdot,\cdot]=0$ & $s(e_3)=e_1$ \\\hline
          ${\bf L^c_{2|2}}$ & $(2|2)$ & $0$ & $X$ odd & $[\cdot,\cdot]=0$ & $s(e_3)=e_1$ \\\hline
          ${\bf L^d_{2|2}}$ & $(2|2)$ & $(\Delta_{2,3},0)$ & $X$ odd & $[e_2,e_3]=X$ & $s(e_3)=e_1$ \\\hline
    \end{tabular}
    }
   \\~\\ Lie superalgebras obtained by extensions of ${\bf L^2_{2|1}}$.
\end{center}
}
\vspace{0.5 em}

{\renewcommand{\arraystretch}{1.5}
\begin{center}
    \makebox[\textwidth][c]{
    \begin{tabular}{|c|c|c|c|c|c|}
    \hline
          \textbf{Name} &\textbf{sdim} & \textbf{Cocycle}& \textbf{Added element} &\textbf{Lie Bracket} &\textbf{Squaring} \\\hline
          ${\bf L^e_{2|2}}$ & $(2|2)$ & $0$ & $X$ even & $[\cdot,\cdot]=0$ & $s(\cdot)=0$ \\\hline
          ${\bf L^e_{1|3}}$ & $(1|3)$ & $0$ & $X$ odd & $[\cdot,\cdot]=0$ & $s(\cdot)=0$ \\\hline
          ${\bf L^f_{2|2}}$ & $(2|2)$ & $(0,\omega_2)$ & $X$ even & $[\cdot,\cdot]=0$ & $s(e_2)=X$ \\\hline
          ${\bf L^g_{2|2}}$ & $(2|2)$ & $(\Delta_{2,3},0)$ & $X$ even & $[e_2,e_3]=X$ & $s(\cdot)=0$ \\\hline
          ${\bf L^f_{1|3}}$ & $(1|3)$ & $(\Delta_{1,2},0)$ & $X$ odd & $[e_1,e_2]=X$ & $s(\cdot)=0$ \\\hline
    \end{tabular}
    }
   \\~\\ Lie superalgebras obtained by extensions of ${\bf L^1_{1|2}}$.
\end{center}
}
\vspace{0.5 em}

{\renewcommand{\arraystretch}{1.5}
\begin{center}
    \makebox[\textwidth][c]{
    \begin{tabular}{|c|c|c|c|c|c|}
    \hline
          \textbf{Name} &\textbf{sdim} & \textbf{Cocycle}& \textbf{Added element} &\textbf{Lie Bracket} &\textbf{Squaring} \\\hline
          ${\bf L^h_{2|2}}$ & $(2|2)$ & $0$ & $X$ even & $[\cdot,\cdot]=0$ & $s(e_2)=e_1$ \\\hline
          ${\bf L^g_{1|3}}$ & $(1|3)$ & $0$ & $X$ odd & $[\cdot,\cdot]=0$ & $s(e_2)=e_1$ \\\hline
          ${\bf L^i_{2|2}}$ & $(2|2)$ & $(0,\omega_3)$ & $X$ even & $[\cdot,\cdot]=0$ & $s(e_2)=e_1, ~s(e_3)=X$ \\\hline
          ${\bf L^j_{2|2}}$ & $(2|2)$ & $(\Delta_{2,3},0)$ & $X$ even & $[e_2,e_3]=X$ & $s(e_2)=e_1$ \\\hline
          ${\bf L^k_{2|2}}$ & $(2|2)$ & $(\Delta_{2,3},\omega_3)$ & $X$ even & $[e_2,e_3]=X$ & $s(e_2)=e_1, ~s(e_3)=X$ \\\hline
          
    \end{tabular}
    }
   \\~\\ Lie superalgebras obtained by extensions of ${\bf L^2_{1|2}}$.
\end{center}
}
\vspace{0.5 em}

{\renewcommand{\arraystretch}{1.5}
\begin{center}
    \makebox[\textwidth][c]{
    \begin{tabular}{|c|c|c|c|c|c|}
    \hline
          \textbf{Name} &\textbf{sdim} & \textbf{Cocycle}& \textbf{Added element} &\textbf{Lie Bracket} &\textbf{Squaring} \\\hline
          ${\bf L^l_{2|2}}$ & $(2|2)$ & $0$ & $X$ even & $[e_2,e_3]=e_1$ & $s(\cdot)=0$ \\\hline
          ${\bf L^h_{1|3}}$ & $(1|3)$ & $0$ & $X$ odd & $[e_2,e_3]=e_1$ & $s(\cdot)=0$ \\\hline
          ${\bf L^m_{2|2}}$ & $(2|2)$ & $(0,\omega_2)$ & $X$ even & $[e_2,e_3]=e_1$ & $s(e_2)=X$ \\\hline
          ${\bf L^n_{2|2}}$ & $(2|2)$ & $(0,\omega_2)+(0,\omega_3)$ & $X$ even & $[e_2,e_3]=e_1$ & $s(e_2)=X, ~s(e_3)=X$ \\\hline
    \end{tabular}
    }
   \\~\\ Lie superalgebras obtained by extensions of ${\bf L^3_{1|2}}$.
\end{center}
}
\vspace{0.5 em}

\vspace{0.5 em}

{\renewcommand{\arraystretch}{1.5}
\begin{center}
    \makebox[\textwidth][c]{
    \begin{tabular}{|c|c|c|c|c|c|}
    \hline
          \textbf{Name} &\textbf{sdim} & \textbf{Cocycle}& \textbf{Added element} &\textbf{Lie Bracket} &\textbf{Squaring} \\\hline
          ${\bf L^o_{2|2}}$ & $(2|2)$ & $0$ & $X$ even & $[e_1,e_2]=e_3$ & $s(\cdot)=0$ \\\hline
          ${\bf L^i_{1|3}}$ & $(1|3)$ & $0$ & $X$ odd & $[e_1,e_2]=e_3$ & $s(\cdot)=0$ \\\hline
          ${\bf L^p_{2|2}}$ & $(2|2)$ & $(0,\omega_2)$ & $X$ even & $[e_1,e_2]=e_3$ & $s(e_2)=X$ \\\hline
          ${\bf L^q_{2|2}}$ & $(2|2)$ & $(0,\omega_3)$ & $X$ even & $[e_1,e_2]=e_3$ & $s(e_3)=X$ \\\hline
          ${\bf L^j_{1|3}}$ & $(1|3)$ & $(\Delta_{1,3},0)$ & $X$ odd & $[e_1,e_2]=e_3, ~[e_1,e_3]=X$ & $s(\cdot)=0$ \\\hline
    \end{tabular}
    }
   \\~\\ Lie superalgebras obtained by extensions of ${\bf L^4_{1|2}}$.
\end{center}
}
\normalsize

\vspace{0.5em}

\subsection{Removing Isomorphic Cases.} As in Section \ref{iso3}, we now remove the isomorphic pairs of algebras from the tables in Section \ref{5.2}. We identify first the pairs which are isomorphic up to a simple relabeling and rearranging of the basis elements. 

\vspace{0.5 em}

\underline{$\sdim(L)=(4|0)$}: $\bf L^b_{4|0} \cong L^c_{4|0}$.\\

\underline{$\sdim(L)=(3|1)$}: $\bf L_{3|1}^a \cong L_{3|1}^c$, $\bf L^b_{1|3} \cong L^d_{1|3}$, $\bf L^e_{1|3} \cong L^g_{1|3}$.\\

\underline{$\sdim(L)=(2|2)$}: $\bf L^a_{2|2}\cong L^e_{2|2}$, $\bf L^b_{2|2} \cong L^o_{2|2}$, $\bf L^c_{2|2} \cong L^f_{2|2} \cong L^h_{2|2}$, $\bf L^d_{2|2} \cong L^p_{2|2}$, $\bf L^g_{2|2} \cong L^l_{2|2}$, 

$\bf L^j_{2|2} \cong L^m_{2|2}$.\\

\underline{$\sdim(L)=(1|3)$}: $\bf L^a_{1|3} \cong L^e_{1|3}$, $\bf L^b_{1|3} \cong L^h_{1|3}$, $\bf L^d_{1|3} \cong L^g_{1|3}$, $\bf L^f_{1|3} \cong L^i_{1|3}$.

\vspace{0.5 em}

\subsubsection{Further Computations}

Now for each superdimension we use invariants such as the dimension of the derived superalgebras and superdimension of the cohomology groups, to deduce a conclusive, non-redundant list of nilpotent Lie superalgebras in dimension 4.

\vspace{0.5 em}

\underline{$\bf L_{4|0}$}: There is only one Lie superalgebra in this case, hence, it is the sole nilpotent one.

\vspace{0.5 em}

\underline{$\bf L_{4|0}$}: Clearly, we find that all three of the remaining Lie superalgebras are non-isomorphic owing to differing dimensions of their derived superalgebras: $\langle[L,L] \cup s(L_{\od})\rangle$. 

\vspace{0.5 em}

\underline{$\bf L_{1|3}$}: Clearly $\bf L^a_{1|3}$ is non-isomorphic to the rest owing to the fact that it is the only superalgebra in this case to have a derived superalgebra of dimension 0. As for the rest, we compute the following table using \texttt{SuperLie} (see \cite{Kr}): 

\vspace{0.5 em}

{\renewcommand{\arraystretch}{1.5}
\footnotesize
\begin{center}
    \begin{tabular}{|c|c|c|c|c|c|}
    \hline
        $L$ & ${\bf L^b_{1|3}}$ & ${\bf L^c_{1|3}}$ & ${\bf L^d_{1|3}}$ & ${\bf L^f_{1|3}}$ & ${\bf L^j_{1|3}}$ \\\hline
        $\sdim\left(H^2(L;\K)\right)$ & $5|0$ & $5|0$ & $5|0$ & $4|2$ & $3|1$ \\\hline
    \end{tabular}
    \\[1em]
    $\sdim\left(H^2(L;\K)\right)$ for Lie superalgebras of $\sdim=(1|3)$.
\end{center}}
\normalsize
\vspace{0.2 em}

It remains to show the relationship between $\bf L^b$, $\bf L^c$ and $\bf L^d$. We find explicit contradictions yielded by the conditions imposed by an isomorphism map $\phi$ to conclude that all three are non-isomorphic to each other. In each case, one can verify that it is not possible to fulfil the conditions without implying $\det(\phi)=0$. 

\vspace{0.5 em}

\underline{$\bf L_{3|1}$}: Again, as in the previous case $\bf L^a$ is non-isomorphic to the other superalgebras. We distinguish between $\bf L^b$, $\bf L^e$ and $\bf L^f$. We find rather trivially through an explicit contradiction that $\bf L^e$ cannot be isomorphic to the other two cases. In particular, the lack of a non-zero bracket imposes conditions that necessitate $\det(\phi)=0$. As for the remaining two cases, we find that $\dim\left(H^2(\mathbf{L^b};\K)\right)=5$, whereas $\dim\left(H^2(\mathbf{L^f};\K)\right)=3$, and hence they are also non-isomorphic.

\vspace{0.5 em}

\underline{$\bf L_{2|2}$}: As with the previous two cases, we find that $\bf L^a$ is non-isomorphic to the other superalgebras. We compute the following table for the superdimension of the second cohomology spaces: 

{\footnotesize
\begin{center}
    \begin{tabular}{|c|c|c|c|c|c|c|c|c|c|}
    \hline
        $L$ & ${\bf L^b_{2|2}}$ & ${\bf L^c_{2|2}}$ & ${\bf L^d_{2|2}}$ & ${\bf L^g_{2|2}}$ & ${\bf L^i_{2|2}}$ & ${\bf L^j_{2|2}}$ & ${\bf L^k_{2|2}}$ & ${\bf L^n_{2|2}}$ & ${\bf L^q_{2|2}}$ \\\hline
        $\sdim\left(H^2(L;\K)\right)$ & $3|2$ & $2|2$ & $2|1$ & $2|2$ & $1|0$ & $1|1$ & $1|0$ & $1|0$ & $1|1$ \\\hline
    \end{tabular}
    \\[1em]
    $\sdim\left(H^2(L;\K)\right)$ for Lie superalgebras of $\sdim=(2|2)$.
\end{center}}
\vspace{0.2 em}

Indeed, using this table we group algebras with the same superdimension as any potential isomorphisms may only arise amongst such pairs. One may derive explicit contradictions for the following cases: $\bf L^c \ncong L^g, ~L^i \ncong L^k, ~L^j \ncong L^q$. Finally, we find that $\bf L^n \cong L^k$, with an isomorphism $\bf \phi: L^n_{2|2} \longrightarrow L^k_{2|2}$ given by:
$$\phi(e_1)=X,~\phi(X)=e_1,~\phi(e_2)=e_2,~\phi(e_3)=e_2+e_3.$$

\vspace{0.5 em}

\noindent We are now ready to list the full classification of dimension 4 nilpotent Lie superalgebras.

\sssbegin{Theorem}\label{classif4} 
The classification of $4$-dimensional nilpotent Lie superalgebras over an algebraically closed field of characteristic $2$ is given by:

\vspace{0.5 em}

\noindent\underline{$\sdim(L)=(0|4)$}:

$\mathbf{L^{1}_{0|4}} = \langle 0|e_{1},e_{2},e_3,e_4\rangle \ (abelian)$.

\vspace{0.5 em}

\noindent\underline{$\sdim(L)=(1|3)$}:

$\mathbf{L^{1}_{1|3}} = \langle e_1|e_{2},e_{3},e_4\rangle \ (abelian)$

$\mathbf{L^{2}_{1|3}} = \langle e_1|e_{2},e_{3},e_4;[e_2,e_3]=e_1\rangle$

$\mathbf{L^{3}_{1|3}} = \langle e_1|e_{2},e_{3},e_4;[e_2,e_3]=e_1, s(e_4)=e_1\rangle$

$\mathbf{L^{4}_{1|3}} = \langle e_1|e_{2},e_{3},e_4; s(e_2)=e_1\rangle$

$\mathbf{L^{5}_{1|3}} = \langle e_{1}|e_{2},e_{3},e_4; [e_1,e_2]=e_4\rangle$

$\mathbf{L^{6}_{1|3}} = \langle e_{1}|e_{2},e_3,e_4; [e_1,e_2]=e_3, [e_1,e_3]=e_4\rangle$

\vspace{0.5 em}

\noindent\underline{$\sdim(L)=(2|2)$}:

$\mathbf{L^{1}_{2|2}} = \langle e_{1},e_{2}|e_3,e_4\rangle \ (abelian)$

            $\mathbf{L^{2}_{2|2}} = \langle e_{1},e_{2}|e_3,e_4; [e_1,e_3]=e_4\rangle$

            $\mathbf{L^{3}_{2|2}} = \langle e_{1},e_{2}|e_{3},e_4; s(e_{3}) = e_{1}\rangle$

            $\mathbf{L^{4}_{2|2}} = \langle e_{1},e_{2}|e_{3},e_4; s(e_{3}) = e_{1}, [e_2,e_3]=e_4\rangle$

            $\mathbf{L^{5}_{2|2}} = \langle e_{1},e_2|e_{3},e_{4}; [e_3,e_4]=e_2\rangle$

            $\mathbf{L^{6}_{2|2}} = \langle e_1,e_2|e_{3},e_{4}; s(e_3)=e_1, s(e_4)=e_2 \rangle$

            $\mathbf{L^{7}_{2|2}} = \langle e_1,e_2|e_{3},e_{4}; s(e_3)=e_1, [e_3,e_4]=e_2 \rangle$

            $\mathbf{L^{8}_{2|2}} = \langle e_1,e_2|e_{3},e_{4}; s(e_3)=e_1, s(e_4)=e_2, [e_3,e_4]=e_2 \rangle$

            $\mathbf{L^{9}_{2|2}} = \langle e_{1},e_2|e_{3},e_4; [e_1,e_3]=e_4, s(e_4)=e_2\rangle$

\vspace{0.5 em}

\noindent\underline{$\sdim(L)=(3|1)$}:

$\mathbf{L^{1}_{3|1}} = \langle e_{1},e_{2},e_3|e_4\rangle \ (abelian)$

            $\mathbf{L^{2}_{3|1}} = \langle e_{1},e_{2},e_3|e_4; [e_1,e_2]=e_3\rangle$

            $\mathbf{L^{3}_{3|1}} = \langle e_{1},e_{2},e_3|e_{4}; s(e_{4}) = e_{1}\rangle$

            $\mathbf{L^{4}_{3|1}} = \langle e_{1},e_{2},e_{3}|e_{4}; [e_1,e_2]=e_3, s(e_{4}) = e_{3}\rangle$

\vspace{0.5 em}

\noindent\underline{$\sdim(L)=(4|0)$}:

$\mathbf{L^{1}_{4|0}} = \langle e_{1},e_{2},e_3,e_4|0\rangle \ (abelian)$

            $\mathbf{L^{2}_{4|0}} = \langle e_{1},e_{2},e_3,e_4|0; [e_1,e_2]=e_4\rangle$

            $\mathbf{L^{3}_{4|0}} = \langle e_{1},e_{2},e_3,e_4|0; [e_1,e_2]=e_3, [e_1,e_3]=e_4\rangle$

\end{Theorem}

\section*{Acknowledgments} It is rare for anyone to accomplish anything worthwhile alone, and this paper is no exception. Growing out of my Senior Year Thesis at New York University Abu Dhabi, the work owes much to the generosity of two remarkable people.

I would like to thank Professor Sofiane Bouarroudj for introducing me to the research question, and more generally the field (pun intended) as a whole. Even more so, I would like to thank him for the faith and encouragement he has shown in me as a supervisor throughout this process. Moreover, I would like to express my immense gratitude to Doctor Quentin Ehret, who has borne the brunt of my incessant questions and discussions, all the while doing so with an alacrity that I will warmly recall.


\end{document}